\pdfoutput=1

\documentclass[11pt]{article}

\usepackage[natbibapa,nodoi]{apacite}
\usepackage{amsmath,amssymb,amsfonts,amsbsy,amsthm,graphicx}

\usepackage[left=1in,right=1in,top=1in,bottom=1in]{geometry}

\theoremstyle{plain}
\newtheorem{theorem}{Theorem}[section]
\newtheorem{lemma}[theorem]{Lemma}

\newtheorem{proposition}[theorem]{Proposition}
\newtheorem{assumption}[theorem]{Assumption}

\theoremstyle{definition}

\theoremstyle{remark}
\newtheorem{remark}{Remark}

\usepackage{dsfont}
\usepackage{algorithm, algorithmic, multirow}
\usepackage{color}

\newcommand{\edit}[1]{{#1}}
\newcommand{\E}{\mathbb{E}}

\begin{document}

\title{Solving Bayesian Risk Optimization via Nested \\ Stochastic Gradient Estimation}

\author{
Sait Cakmak, Di Wu, and Enlu Zhou\\
School of Industrial and Systems Engineering \\
Georgia Institute of Technology\\
Atlanta, GA 30332
}

\maketitle

\begin{abstract}
    In this paper, we aim to solve Bayesian Risk Optimization (BRO), which is a recently proposed framework that formulates simulation optimization under input uncertainty.
    In order to efficiently solve the BRO problem, we derive nested stochastic gradient estimators and propose corresponding stochastic approximation algorithms. 
    We show that our gradient estimators are asymptotically unbiased and consistent, and that the algorithms converge asymptotically. 
    We demonstrate the empirical performance of the algorithms on a two-sided market model. 
    Our estimators are of independent interest in extending the literature of stochastic gradient estimation to the case of nested risk functions.
\end{abstract}


\section{Introduction}
We consider the following optimization problem:
\begin{equation} \label{sim-opt-problem1}
    \min_{x \in \mathcal{X}} H(x) := \E_{\xi \sim \mathbb{P}^\text{c}} [h(x, \xi)],
\end{equation}
where $\mathcal{X}$ is the solution space, $\xi$ is a random vector representing the randomness in simulation, and $h(\cdot, \cdot)$ is a function that is evaluated through simulation. The expectation is taken with respect to (w.r.t.) $\mathbb{P}^\text{c}$, the correct distribution of $\xi$. In a typical simulation optimization setting, the true distribution $\mathbb{P}^\text{c}$ is unknown and estimated from a finite set of input data, and the following approximate problem is solved, where the estimated distribution is denoted by $\hat{\mathbb{P}}$. 
\begin{equation} \label{approximate-problem}
    \min_{x \in \mathcal{X}} \E_{\xi \sim \hat{\mathbb{P}}} [h(x, \xi)]
\end{equation}
Due to the use of a finite dataset, even when the approximate problem (\ref{approximate-problem}) is solved to optimality, the optimal solution can perform poorly under true distribution. This issue is referred to as input uncertainty in simulation optimization.

Recently, \cite{zhou2015BRO} and \cite{wu2018BRO} proposed the Bayesian Risk Optimization (BRO) framework which formulates the simulation optimization problem under input uncertainty. 
In BRO, assuming that $\mathbb{P}^\text{c}$ belongs to a known parameterized family of distributions $\{\mathbb{P}_{\theta}\}_{\theta \in \Theta}$ with unknown parameter $\theta^{\text{c}}$, instead of solving (\ref{sim-opt-problem1}) we solve the following:
\begin{equation} \label{bro-problem}
    \min_{x \in \mathcal{X}}  \rho_{\theta \sim \mathbb{P}^N} \{H(x; \theta)\}= \rho_{\theta} \{ \E_{\xi \sim \mathbb{P}_{\theta}}[h(x, \xi)] \},
\end{equation}
where $\rho$ is a risk function mapping the random variable $H(x; \theta)$ (induced by $\theta \sim \mathbb{P}^N$) to a real number, and $\mathbb{P}^N$ is the Bayesian posterior distribution of $\theta$ given a chosen prior and input data $\phi^N = \{\zeta_i\}_{i=1}^N$. The risk function $\rho$ can be chosen according to the risk preferences of the practitioner, and includes the risk neutral expectation and the minimax formulation of distributionally robust optimization (DRO, see e.g. \citealp{rahimian2019dro-review}) as extreme cases under certain conditions. In this paper, we consider the following four cases of $\rho$: Expectation, Mean-Variance, Value-at-Risk (VaR), and Conditional Value-at-risk (CVaR). A formal introduction and a thorough review of BRO, along with a discussion on alternative approaches, is provided in Section \ref{bro-review}.

We aim to solve the BRO problem (\ref{bro-problem}). To do so, we will use a Stochastic Approximation (SA, see \citealp{Kushner&Yin2003SA}) approach, which requires gradient information. Historically, most work on stochastic gradient estimation focused on finding the gradient of expectation (see \citealp{FU2006Gradient, FU2008SA-Summary}). Some more recent research studies the Monte-Carlo estimation of gradients of VaR and CVaR; e.g. \cite{Hong2009VaR} for VaR, and \cite{Hong2009CVaR} for CVaR, where each derives a closed form expression of the corresponding gradient, and provides an asymptotically unbiased and consistent infinitesimal perturbation analysis (IPA) based estimator. Other work in this field includes  \cite{Liu&Hong2009Quantile-Kernel}, \cite{Fu2009Quantile-Sens}, \cite{Jiang2015Quantile}, \cite{Tamar2015CVaR} and \cite{Peng2017Quantile}, to name a few. A review of Monte-Carlo methods for estimation of VaR, CVaR and their gradients can be found in \cite{Hong2014VaR-CVaR-Rev}.

Other related works include the literature on nested simulation, e.g. \cite{Lan2010CVaR}, \cite{Gordy2010Nested}, \cite{Broadie2015Regression}, \cite{Zhu2017input-uncertainty}; \cite{Jaiswal2019}, which studies the data-driven risk averse optimization problem under a parameterized Bayesian setting using the log-exponential risk measure; and \cite{Wang2020BOBRO}, which uses Bayesian Optimization (see \citealp{frazier2018tutorial}) methods to optimize the expectation case of BRO with black-box expensive-to-evaluate objective functions.

Our work differs from the aforementioned works in the sense that the literature on nested simulation does not consider gradients or optimization, and the literature on gradient estimation does not consider nested risk functions. The literature on gradient estimation requires access to $H(x; \theta)$ (and its gradients), while we only have access to $h(x, \xi(\theta))$ (and its gradients), where $H(x; \theta) := \E_{\xi \sim \mathbb{P}_{\theta}} [h(x, \xi)]$, and $H(x; \theta)$ (and its gradients) has to be estimated via sampling. The need to estimate the function $H(x; \theta)$ adds another level of uncertainty to gradient estimation. To the best of our knowledge, this work is the first to study stochastic gradient estimation of nested risk functions. 

The contributions of this paper can be summarized as follows: (i) We propose sample path gradient estimators $\rho \{H(x; \theta) \}$ for the four risk functions $\rho$ mentioned earlier, extending the literature in stochastic gradient estimation to the case of nested risk functions; (ii) We propose stochastic approximation algorithms with local convergence guarantees for optimization of the BRO framework; (iii) We provide a numerical study on a two-sided market model that demonstrates the value of risk averse solution approaches in the presence of input uncertainty. Although the exposition in this paper is focused on the BRO framework, it is worth noting that our estimators can be applied more broadly, e.g., for estimating the sensitivities of quantiles of financial portfolios.

\section{An Overview of BRO Framework} \label{bro-review}
As mentioned in the introduction, in a typical simulation optimization framework, one aims to solve the following problem:
\begin{equation} \label{simopt-prob-2}
    \min_{x \in \mathcal{X}} H(x) := \E_{\xi \sim \mathbb{P}^\text{c}}[h(x, \xi)],
\end{equation}
where the solution space $\mathcal{X}$ is a non-empty, compact subset of $\mathbb{R}^{d_1}$, $\xi \in \mathbb{R}^{d_2}$ is a random vector representing the stochastic noise in the system, and $h$ is a function mapping $\mathbb{R}^{d_1} \times \mathbb{R}^{d^2}$ to $\mathbb{R}$. The expectation is taken w.r.t $\mathbb{P}^\text{c}$, the true distribution of $\xi$. In practice, $\mathbb{P}^\text{c}$ is not known and is typically replaced with an estimate $\hat{\mathbb{P}}$ which is obtained from a finite set of input data. 
The estimation error of $\hat{\mathbb{P}}$ due to the use of finite data is often referred to as the input model uncertainty, or simply as input uncertainty. There is a large literature dedicated to studying the impact of input uncertainty in estimating system performance; see \cite{Barton2012Input-uncertainty} and \cite{Song2014Input-uncertainty} for a review.

Due to input uncertainty, even when the estimated problem $\min_{x} \E_{\xi \sim \hat{\mathbb{P}}} [h(x, \xi)]$ is solved to optimality, the optimal solution can perform poorly under the true distribution. Hence, a natural question is, ``\textit{how do we make decisions that account for input uncertainty?}". We are interested in finding good solutions that hedge against input uncertainty. 
One common approach is to construct an ambiguity set $\mathcal{D}$ that includes $\mathbb{P}^\text{c}$ with high probability, and optimize w.r.t. the worst-case outcome within this set. This approach is referred to as Distributionally Robust Optimization (DRO) framework and has a large literature dedicated to it; see \cite{rahimian2019dro-review} for a review. In DRO, constructing the ambiguity set is a non-trivial task and has a large impact on the solution performance and tractability of the resulting problem. A large ambiguity set can lead to overly conservative solutions, whereas, a small uncertainty set might fail to include the true distribution. An alternative approach is to optimize with respect to a risk neutral expectation of the objective function over the set of all possible input distributions. 
\edit{
As argued in \cite{zhou2015BRO}, these two approaches can be seen as two extreme cases. The risk neutral expectation might fail to put enough weight over extreme (tail) scenarios, whereas, the DRO approach might be overly conservative due to hedging against worst-case scenarios. 

Moreover, to the best of our knowledge, in a simulation optimization setting, the existing literature lacks tractable reformulations of the DRO problems. Although one can use the minimax formulation to formulate the distributionally robust simulation optimization problem, efficient optimization of this problem remains an open question due to the lack of structure in the $h(\cdot,\cdot)$ function. In the robust optimization literature, \cite{Bertsimas2010RobustSim} study a simulation optimization problem that is jointly robust to both implementation errors and parameter uncertainty, however, their method does not work when one is only concerned about the parameter uncertainty. 
Besides these popular approaches, the robust simulation optimization problem found interest in the kriging literature, e.g. \cite{Dellino2015Metamodel, KLEIJNEN2017Kriging}, where a response surface is fitted over $\mathcal{X} \times \Theta$, and robustness is typically facilitated by optimizing the mean performance subject to constraints on the standard deviation. 
}

\edit{
In this paper, we focus on the Bayesian Risk Optimization (BRO) framework, which was proposed by \cite{zhou2015BRO} and \cite{wu2018BRO} as an alternative approach to simulation optimization under input parameter uncertainty.
}
Suppose that the true distribution $\mathbb{P}^\text{c}$ belongs to a parameterized family of distributions $\{\mathbb{P}_{\theta}\}_{\theta \in \Theta}$ such that $\mathbb{P}^\text{c} = \mathbb{P}_{\theta^{\text{c}}}$ for some $\theta^{\text{c}}$, where $\theta^{\text{c}} \in \Theta$ is the unknown true parameter and $\Theta$ is the parameter space. Assuming that the form of $\mathbb{P}_{\theta}$ is known, we take a Bayesian approach and calculate the posterior likelihood of $\theta$ for a given dataset $\phi^N := \{\xi_i\}_{i=1}^N$ of independent and identically distributed (i.i.d.) input data drawn from the true distribution. Let $p(\theta)$ denote the prior distribution of $\theta^{\text{c}}$. Then, using Bayesian updating we can calculate the posterior distribution
\begin{equation*}
    \mathbb{P}^N := p(\theta \mid \phi^N) \propto p(\theta)p(\phi^N \mid \theta) = p(\theta) \prod_{i=1}^N f(\xi_i \mid \theta),
\end{equation*}
where $p(\phi^N \mid \theta)$ ($f(\xi_i \mid \theta)$) is the likelihood of obtaining $\phi^N$ ($\xi_i$) given parameter $\theta$, and $\propto$ denotes equivalence up to a normalization constant. 

Define $H(x; \theta) := \E_{\mathbb{P}_{\theta}} [h(x, \xi)]$ as the objective value under parameter $\theta$, where $\E_{\mathbb{P}_{\theta}}$ denotes the expectation w.r.t. $\xi \sim \mathbb{P}_{\theta}$. If we view $\theta$ as a random variable with distribution $\mathbb{P}^N$, we can treat $H(x; \theta)$ as a random variable induced by $\theta$. Define $\rho$ as a risk function over $H(x; \theta)$ which maps the random variable to $\mathbb{R}$. Instead of solving (\ref{simopt-prob-2}), we solve
\begin{equation} \label{bro-prob2}
    \min_{x \in \mathcal{X}}  \rho_{\theta \sim \mathbb{P}^N} \{H(x; \theta)\} = \rho_{\theta} \{ \E_{\mathbb{P}_{\theta}}[h(x, \xi)] \},
\end{equation}
which is referred to as the BRO problem. The risk function $\rho$ can be chosen to reflect the risk preference of the practitioner. In this paper, we focus on the following four cases of $\rho$:
\begin{enumerate}
    \item Expectation: $ \min_{x \in \mathcal{X}} \E_{\theta} \left[ \E_{\mathbb{P}_{\theta}} [h(x,\xi)] \right]$;
    
    \item Mean - Variance: $ \min_{x \in \mathcal{X}} \E_{\theta} \left[ \E_{\mathbb{P}_{\theta}} [h(x,\xi)] \right] + a \text{Var}_{\theta} \left( \E_{\mathbb{P}_{\theta}} [h(x,\xi)] \right) $;
    
    \item Value-at-Risk: $ \min_{x \in \mathcal{X}} \text{VaR}_{\alpha} \left( \E_{\mathbb{P}_{\theta}} [h(x,\xi)] \right) $;
    
    \item Conditional Value-at-Risk: $ \min_{x \in \mathcal{X}} \text{CVaR}_{\alpha} \left( \E_{\mathbb{P}_{\theta}} [h(x,\xi)] \right) $;
\end{enumerate}
where $\E_{\theta}$ ($\text{Var}_{\theta}$) denote that the expectation (variance) is taken w.r.t. $\theta \sim \mathbb{P}^N$, $\text{VaR}_{\alpha}$ and $\text{CVaR}_{\alpha}$ denote the $\alpha$ level Value-at-Risk and Conditional Value-at-Risk respectively. We will define VaR and CVaR formally in corresponding subsections.

For the four cases of $\rho$ considered here, \cite{wu2018BRO} study the asymptotic properties of the objective functions and optimal solutions. We briefly summarize their results here. As the intuition would suggest, they show that as the data size $N \rightarrow \infty$, the posterior distribution $\mathbb{P}^N$ converges in distribution to a degenerate distribution on $\theta^{\text{c}}$. Furthermore, under mild regularity conditions, it is shown that for every fixed $x \in \mathcal{X}$ as $N \rightarrow \infty$, $\rho_{\theta} \{H(x; \theta)\} \rightarrow H(x; \theta^{\text{c}})$ almost surely (a.s.), and $\min_{x \in \mathcal{X}} \rho_{\theta} \{H(x; \theta)\} \rightarrow \min_{x \in \mathcal{X}} H(x; \theta^{\text{c}})$ a.s. for all four cases of $\rho$ considered here. Similarly, for the consistency of optimal solutions, it is shown that $\mathbb{D}(S_N, S) \rightarrow 0$ a.s. as $N \rightarrow \infty$ where $S_N$ and $S$ are the sets of optimal solutions to (\ref{bro-prob2}) and (\ref{simopt-prob-2}) respectively, and $\mathbb{D}(A, B) := \sup_{x \in A} dist(x, B)$ is the distance between two sets with $dist(x, B) := \inf_{y \in B} \|x - y\|$ and $\|.\|$ being an arbitrary norm.

Moreover, the analysis of \cite{wu2018BRO} reveals the following asymptotic normality results which can be used to construct confidence intervals for the true objective value. Let $\mathcal{N}$ denote the normal distribution, and let $\phi$ and $\Phi$ denote the probability density function (PDF) and cumulative density function (CDF) of $\mathcal{N}(0,1)$ respectively. Then, for every $x \in \mathcal{X}$, as $N \rightarrow \infty$,
\begin{itemize}
    \item for Expectation and Mean-Variance objectives,
    \begin{equation*}
        \sqrt{N} \{ \E_{\theta}[H(x; \theta)] + a \text{Var}[H(x; \theta)] - H(x; \theta^{\text{c}}) \} \Rightarrow \mathcal{N}(0, \sigma_{x}^2);
    \end{equation*}
    \item for the Value-at-Risk objective,
    \begin{equation*}
        \sqrt{N} \{ \text{VaR}_{\alpha}[H(x; \theta)] - H(x; \theta^{\text{c}}) \} \Rightarrow \mathcal{N}(\sigma_{x} \Phi^{-1}(\alpha), \sigma_{x}^2);
    \end{equation*}
    \item and for the Conditional Value-at-Risk objective,
    \begin{equation*}
        \sqrt{N} \{ \text{CVaR}_{\alpha}[H(x; \theta)] - H(x; \theta^{\text{c}}) \} \Rightarrow \mathcal{N} \left(\frac{\sigma_{x}}{1 - \alpha} \phi(\Phi^{-1}(\alpha)), \sigma_{x}^2 \right);
    \end{equation*}
\end{itemize}
where $\Rightarrow$ denotes convergence in distribution. The variance is defined as \linebreak $\sigma_{x}^2 := \nabla_{\theta} H(x; \theta^{\text{c}})^{\top} [\mathcal{I}(\theta^{\text{c}})]^{-1} \nabla_{\theta} H(x; \theta^{\text{c}})$ where $\mathcal{I}(\theta^{\text{c}})$ is the Fisher information matrix, $\top$ denotes the transpose, and $\nabla_{\theta}$ is the gradient w.r.t. $\theta$. The point-wise convergence results presented here can be extended to convergence results in the function space of $H(\cdot; \theta)$. Similar normality results also hold for the optimal values.

To summarize, \cite{wu2018BRO} establish the consistency and asymptotic normality of objective functions and optimal solutions for the four cases of $\rho$ considered here. They also show that the objectives of BRO can be approximated as a weighted sum of posterior mean objective and half-width of the true-objective's confidence interval. 
In this paper, our aim is to {\it optimize the BRO problem} (\ref{bro-prob2}) for a given choice of $\rho$ and a given posterior distribution $\mathbb{P}^N$ of $\theta$. 
\edit{
We refer the interested reader to \cite{zhou2015BRO}, \cite{Zhou2017book-chapter}, and \cite{wu2018BRO} for further discussion on BRO formulation.
}

\section{Solving the BRO Problem} \label{estimator-section}
In this section, we introduce our approach to solving the BRO problem. We take an SA approach, develop the stochastic gradient estimators needed, and conclude with convergence results for the algorithms. 
\edit{Throughout the paper, we use $d(\cdot, \cdot)$ and $D(\cdot, \cdot)$ to denote the gradients $\frac{d h(\cdot, \cdot)}{d x}$ and $\frac{d H(\cdot, \cdot)}{dx}$ respectively. 
We use $\E_{\theta}$ as shorthand for $\E_{\theta \sim \mathbb{P}^N}$, the expectation over the posterior distribution of $\theta$, and $\E_{\mathbb{P}_{\theta}}$ as a shorthand for $\E_{\xi \sim \mathbb{P}_{\theta}}$. 
The nested expectation $\E_{\theta}[\E_{\mathbb{P}_{\theta}}[\cdot]]$ is also shortened as $\E_{\theta, \mathbb{P}_{\theta}}[\cdot]$.
} 
The proofs are provided in the online supplement.

\subsection{Stochastic Approximation Algorithm}
The BRO problem in its essence is a typical simulation optimization problem where the objective function is costly to estimate. Due to the nested structure of the objective function, one needs many more samples to estimate the BRO objective compared to a typical expectation or CVaR objective. If one were to use $m$ samples to estimate the inner expectation and $n$ samples to estimate the outer risk function, it would take a total of $n \times m$ samples to estimate the BRO objective. This high cost of estimation motivates us to concentrate on algorithms that take advantage of the structure of the problem and require fewer function evaluations per iteration. With this motivation, the class of gradient based methods known as Stochastic Approximation emerges as an obvious candidate. To solve the BRO problem (\ref{bro-prob2}), we propose to use the SA algorithm of the following form (see \citealp{Kushner&Yin2003SA}): 
\begin{equation} \label{rm-algorithm}
x_{t+1} = \Pi_{\mathcal{X}} [x_t + \epsilon_t Y_t]
\end{equation}
where $\mathcal{X}$ is a non-empty, compact solution space, $\{\epsilon_t\}_{t \geq 0}$ is the step size sequence, $Y_t$ is the descent direction, $\Pi$ is the projection operator that projects the iterate back to the feasible set $\mathcal{X}$. A typical candidate for $Y_t$ is an estimate of the negative gradient of the objective function, which leads to the well known Stochastic Gradient Descent algorithm. Given a good estimator of the gradient, the stochastic approximation algorithm has nice convergence properties. We proceed in next subsection with the derivation of stochastic gradient estimators of the BRO problem (\ref{bro-prob2}) for the four cases of $\rho$ mentioned above.

\subsection{Derivation of Stochastic Gradient Estimators}
In this section, we derive the stochastic gradient estimators for the BRO problem. The results are derived only for one-dimensional $x$. Multidimensional case can be handled by treating each dimension as a one-dimensional parameter while fixing the rest. We start by providing the estimators for Expectation and Mean - Variance cases without going into details, then derive the estimators for the more technically challenging cases of VaR and CVaR. The following lemma from \cite{Broadie1996IPA} is key to the consistency of IPA estimators and is used without mention throughout the paper.

\begin{lemma} Proposition 1, \cite{Broadie1996IPA} -
    Let $\phi$ denote a Lipschitz continuous function that is differentiable on a set of points $D_{\phi}$. Suppose that there exists a random variable $K(\xi)$ with $\E[K(\xi)] < \infty$ such that $|h(x_1, \xi) - h(x_2, \xi)| < K(\xi) |x_1 - x_2|$ for all $x_1, x_2 \in \mathcal{X}$ and $d(x, \xi)$ exists w.p. (with probability) 1 for all $x \in \mathcal{X}$, with $\mathcal{X}$ an open set. If $P(h(x, \xi) \in D_{\phi}) = 1$ for all $x \in \mathcal{X}$, then $d \E[\phi(h(x, \xi))] / dx = \E[\phi'(h(x, \xi)) d(x, \xi)]$ for all $x \in \mathcal{X}$.
\end{lemma}

\subsubsection{Expectation and Mean-Variance Cases} 
Suppose that the interchange of gradient and expectation is justified (see Assumption \ref{interchange-assumption}). Then, we have the following for the gradients of expectation and variance respectively:
\begin{equation} \label{expectation-gradient}
    \frac{d \E_{\theta} \left[\E_{\mathbb{P}_{\theta}}[h(x, \xi)] \right]}{d x} = \E_{\theta} \left[\frac{d \E_{\mathbb{P}_{\theta}}[h(x, \xi)]}{d x} \right]  = \E_{\theta} \left[ \E_{\mathbb{P}_{\theta}} \left[\frac{d h(x, \xi)}{d x} \right] \right] = \E_{\theta, \mathbb{P}_{\theta}} \left[ d(x, \xi) \right]
\end{equation}
and 
\begin{equation} \label{variance-gradient}
\begin{aligned}
    \frac{d \text{Var}_{\theta} \left(\E_{\mathbb{P}_{\theta}}[h(x, \xi)] \right)}{d x} &= \frac{d \left( \E_{\theta} \left[ \E_{\mathbb{P}_{\theta}} [h(x, \xi)]^2 \right] - \E_{\theta} \left[ \E_{\mathbb{P}_{\theta}} [h(x, \xi)] \right]^2 \right)}{dx} \\
    &= \E_{\theta} \left[ \frac{d \E_{\mathbb{P}_{\theta}} [h(x, \xi)]^2 }{dx} \right] - 2 \E_{\theta} \left[ \E_{\mathbb{P}_{\theta}} [h(x, \xi)] \right] \frac{d \E_{\theta} \left[ \E_{\mathbb{P}_{\theta}} [h(x, \xi)] \right]}{dx} \\
    &= 2 \left( \E_{\theta} \left[\E_{\mathbb{P}_{\theta}}[h(x, \xi)] \E_{\mathbb{P}_{\theta}}[d(x, \xi)] \right] - \E_{\theta, \mathbb{P}_{\theta}}[ h(x, \xi) ] \E_{\theta, \mathbb{P}_{\theta}}[ d(x, \xi) ] \right).
\end{aligned}
\end{equation}
The equations (\ref{expectation-gradient}) and (\ref{variance-gradient}) can be used to provide gradient estimators for Expectation and Mean-Variance cases. For the Expectation case, it is seen that $d(x, \xi(\theta))$ is a single run unbiased gradient estimator and the sample average $\frac{1}{n}\sum_{i=1}^n d(x, \xi_i(\theta_i))$ (where $\xi_i(\theta_i)$ are independent with distribution $\mathbb{P}_{\theta_i}$ with $\theta_i \stackrel{iid}{\sim} \mathbb{P}^N$) is a strongly consistent estimator of the gradient. Similarly, for Mean-Variance case, we have
\begin{equation} \label{mean-variance-estimator}
    d(x, \xi_1(\theta_1)) + 2a \left( h(x, \xi_2(\theta_2)) d(x, \xi_3(\theta_2)) - h(x, \xi_4(\theta_3)) d(x ,\xi_5(\theta_4)) \right)
\end{equation}
as an unbiased gradient estimator with $\xi_1, \ldots, \xi_5$ independent and $\theta_1, \ldots, \theta_4$ i.i.d. samples. One could use the same sample $\theta$ for $\theta_1, \theta_2 \ \& \ \theta_3$ and the same sample $\xi$ for $\xi_1, \xi_2 \ \& \ \xi_4$ at the expense of increased variance, and reduce the number of simulation runs to $3$. However, using any fewer simulation runs would make the estimation of second and third terms biased. We do not study the trade off here since our main focus is on the estimation of VaR and CVaR gradients. This subsection is concluded by noting that sample averaging yields a strongly consistent estimator for the Mean-Variance case.

\subsubsection{Value-at-Risk Case}
In this subsection, we introduce the nested estimator of VaR gradients, and establish the asymptotical properties of the proposed estimator. Value-at-Risk, defined as \linebreak $\text{VaR}_{\alpha} (H(x; \theta)) = \inf \{t: P(H(x; \theta) \leq t) \geq \alpha \} $, is the $\alpha$ quantile of the loss function. We are interested in estimating the gradient $d \text{VaR}_{\alpha} (H(x; \theta))/d x$ using samples of $h(x, \xi(\theta))$ and corresponding sample path gradients. Throughout the paper, $v_{\alpha}(x)$ and $v'_{\alpha}(x)$ are used as shorthand notations for VaR$_{\alpha}(H(x; \theta))$ and $d \text{VaR}_{\alpha}(H(x; \theta)) / dx$ respectively.

If one has access to $n$ samples of $H(x; \theta)$, then $v_{\alpha}(x)$ can be estimated by the sample quantile $\hat{v}^n_{\alpha}(x) := H(x; \theta_{(\lceil \alpha n \rceil)})$ (see \citealp{Serfling1980}) where $\lceil \cdot \rceil$ is the ceiling function, $\theta_{(i)}$ denotes $i^{th}$ order statistic corresponding to the ordering $H(x; \theta_{(1)}) \leq H(x; \theta_{(2)}) \leq \ldots \leq H(x; \theta_{(n)})$, and $H$ is treated as a random variable induced by $\theta \sim \mathbb{P}^N$. However, in our case, we only have access to samples from $h(x, \xi(\theta))$. Let $\hat{H}^m(x; \theta) := \frac{1}{m} \sum_{j=1}^m h(x, \xi_j(\theta))$ denote the Monte-Carlo estimator of $H(x; \theta)$ generated using $m$ samples. 
Note that the ordering of $\theta_{(i)}$ based on $H$ does not necessarily correspond to the ordering of $\hat{H}^m$, i.e. $\hat{H}^m (x; \theta_{(1)}) \leq \hat{H}^m (x; \theta_{(2)}) \leq \ldots \leq \hat{H}^m (x; \theta_{(n)})$ does not hold in general. Therefore, we define a new ordering, denoted by $\hat{\theta}^m_{(i)}$, such that $\hat{H}^m (x; \hat{\theta}^m_{(1)}) \leq \hat{H}^m (x; \hat{\theta}^m_{(2)}) \leq \ldots \leq \hat{H}^m (x; \hat{\theta}^m_{(n)})$. This ordering is not uniquely defined by $\{\theta_1, \theta_2, \ldots, \theta_n\}$ and depends on the realization of $\xi(\theta)$s. Under a mild set of assumptions, \cite{Zhu2017input-uncertainty} shows that $\hat{v}^{n,m}_{\alpha}(x) := \hat{H}^m(x; \hat{\theta}^m_{(\lceil \alpha n \rceil)})$ is a strongly consistent estimator of $\text{VaR}_{\alpha}$. Motivated by the consistency of $\hat{v}^{n,m}_{\alpha}(x)$, we propose 
\begin{equation}
    \varphi^{n,m}_{\alpha}(x) := \partial_x \hat{H}^m (x; \theta) |_{ \hat{H}^m (x; \theta) = \hat{v}^{n,m}_{\alpha} } = \hat{D}^m (x; \hat{\theta}^m_{(\lceil \alpha n \rceil)}) 
\end{equation}
as the nested estimator of VaR gradients where $\hat{D}^m (x; \theta) := \frac{1}{m} \sum_{j=1}^m d(x, \xi_j(\theta))$ is the IPA gradient estimator corresponding to $\hat{H}^m (x; \theta)$. In the remainder of this subsection, we proceed to show that the estimator $\varphi^{n,m}_{\alpha}(x)$ is asymptotically unbiased, and the batch-mean estimator $\bar{\varphi}^{n,m,k}_{\alpha}(x) := \frac{1}{k} \sum_{i=1}^k \varphi^{n,m}_{\alpha, i}(x)$, where $k$ is the number of batches of equal size and $\varphi^{n,m}_{\alpha, i}(x)$ is the estimator corresponding to batch $i$, is consistent and asymptotically normally distributed.

Notice that $\varphi^{n,m}_{\alpha}(x)$ has the same form as $I^n = \partial_x H(x; \theta) |_{H(x; \theta) = \hat{v}^n_{\alpha}} = D(x; \theta_{(\lceil \alpha n \rceil)})$, the single-layer estimator of quantile gradients of \cite{Hong2009VaR}. Both estimators stem from the observation that, under a mild set of assumptions, the quantile gradients can be expressed as $v'_{\alpha}(x) = \E_{\theta}[\partial_x H(x; \theta) \mid H(x; \theta) = v_{\alpha}(x)]$. 

We now introduce the technical conditions that lead to the consistency of these estimators. The assumptions we introduce here can be viewed in three categories. First, we have a set of assumptions due to \cite{Zhu2017input-uncertainty} that are needed to justify consistency of $\hat{v}^{n,m}_{\alpha}(x)$ by providing the necessary smoothness of $\hat{H}^m(x; \theta)$. A second set of assumptions are needed to justify the interchange of gradient and expectation, and thus the validity of IPA gradient estimators. An additional assumption by \cite{Hong2009VaR} is needed to validate the interchange for the case of VaR. A final set of assumptions are needed to mitigate the difficulties arising from conditioning on measure zero sets, and ensure that the pathwise gradient estimator $d(x, \xi(\theta))$ is sufficiently smooth.

Let $ \mathcal{E} (x, \xi(\theta)) = h(x, \xi(\theta)) - H(x; \theta) $ denote the estimation error and $ \bar{\mathcal{E}}^m (x; \theta) = \sqrt[]{m}  \frac{1}{m} \sum_{j=1}^m \mathcal{E} (x, \xi_j(\theta))$ denote the normalized error. Then, $\hat{H}^m (x; \theta) = H(x; \theta) + \frac{1}{\sqrt{m}} \bar{\mathcal{E}}^m (x; \theta)$. Under following set of assumptions, \cite{Zhu2017input-uncertainty} prove that $\hat{v}^{n,m}_{\alpha} := \hat{H}^m(x; \hat{\theta}^m_{(\lceil \alpha n \rceil)})$ is a strongly consistent estimator of $\text{VaR}_{\alpha}$. 

\begin{assumption} \label{zhu-assumption}
\citep{Zhu2017input-uncertainty}
\begin{enumerate}
    \item For all $x \in \mathcal{X}$, the response $h(x, \xi(\theta))$ has finite conditional second moment, i.e., $\tau_{\theta}^2 = \E_{\mathbb{P}_{\theta}}[h(x, \xi)^2] < \infty$ w.p. 1 $(\mathbb{P}^N)$ and $\tau^2 = \E_{\theta, \mathbb{P}_{\theta}}[h(x, \xi)^2] = \int \tau_{\theta}^2 d\mathbb{P}^N < \infty$.
    
    \item The joint density $p_m(h,e)$ of $H(x; \theta)$ and $\bar{\mathcal{E}}^m(x; \theta)$, and its partial gradients $\frac{d}{dh} p_m(h,e)$ and $\frac{d^2}{dh^2}p_m(h,e)$ exist for each $m$, all pairs of $(h,e)$ and for all $x \in \mathcal{X}$.
    
    \item For all $x \in \mathcal{X}$, there exists non-negative functions $g_{0,m}(\cdot), g_{1,m}(\cdot)$ and $g_{2,m}(\cdot)$ such that $p_m(h,e) \leq g_{0,m}(e)$, $|\frac{d}{dh}p_m(h,e)| \leq g_{1,m}(e)$, $|\frac{d^2}{dh^2}p_m(h,e)| \leq g_{2,m}(e)$ for all $(h,e)$. Furthermore, $\sup_m \int |e|^r g_{i,m}(e) de < \infty $ for $i = 0,1,2$, and $0\leq r\leq 4$.
\end{enumerate}
\end{assumption}

The first part of the assumption ensures the validity of Central Limit Theorem (CLT). Second and third parts first appear in \cite{Gordy2010Nested} and provide sufficient smoothness to ensure that the PDF of $\hat{H}^m(\cdot)$ convergences to the PDF of $H(\cdot)$ sufficiently fast. For our purposes, they provide uniform bounds on the moments of the estimation error. See \cite{Gordy2010Nested} for further discussion on these assumptions.

\begin{assumption} \label{interchange-assumption}
There exists a random variable $K(\xi(\theta))$ such that $\E_{\theta, \mathbb{P}_{\theta}}[K(\xi)] < \infty$, and the following holds in a probability 1 $(\mathbb{P}^N)$ subset of $\Theta$.
\begin{enumerate}
    \item $|h (x_2, \xi(\theta)) - h (x_1, \xi(\theta))| \leq K(\xi(\theta)) |x_2 - x_1|$ w.p.1 $(\mathbb{P}_{\theta})$ for all $x_1, x_2 \in \mathcal{X}$.
    \item The sample path gradient $d(x, \xi(\theta))$ exists w.p.1 $(\mathbb{P}_{\theta})$.
\end{enumerate}
\end{assumption}

Assumption \ref{interchange-assumption} ensures that $D(x; \theta)$ exists (w.p.1), $D(x; \theta) = \E_{\mathbb{P}_{\theta}}[d(x, \xi)]$, and that these relations can be extended to $\E_{\theta}[H(x; \theta)]$. Let $F(\cdot; x)$ denote the distribution function of $H(x; \theta)$ and define $g(t; x) := \E_{\theta}[\partial_x H(x; \theta) \mid H(x; \theta) = t]$. We have the following assumption due to \cite{Hong2009VaR}.

\begin{assumption} \label{hong-var-assumptions} \citep{Hong2009VaR}
    For any $x \in \mathcal{X}$, $H(x, \theta)$ has a continuous density $f(t; x)$ in a neighborhood of $t = v_{\alpha}(x)$, and $\partial_{x} F(t; x)$ exists and is continuous w.r.t. both $x$ and $t$ at $t = v_{\alpha}(x)$.
\end{assumption}

This assumption ensures that $H(\cdot)$ is a continuous random variable in a neighborhood of $v_{\alpha}(x)$, and that its gradient exists and is continuous in the same neighborhood. It is shown in \cite{Jiang2015Quantile} that Assumptions \ref{interchange-assumption} \& \ref{hong-var-assumptions} are sufficient to justify the expression $v'_{\alpha}(x) = \E_{\theta}[\partial_x H(x; \theta) \mid H(x; \theta) = v_{\alpha}(x)]$, and that the continuity (in $t$) of $g(t; x)$ follows from these two assumptions. Under these assumptions, \cite{Hong2009VaR} and \cite{Jiang2015Quantile} show that $\E_{\theta}[I_n] \rightarrow v'_{\alpha}(x)$ as $n \rightarrow \infty$, and the batch-mean estimator $\Bar{I}^{n,k} = \frac{1}{k} \sum_{i=1}^k I^n_i$ (with $k$ as the number of batches) is consistent.

We would like to show that $\E_{\theta, \mathbb{P}_{\theta}} [\varphi^{n,m}_{\alpha}(x)] \rightarrow v'_{\alpha}(x)$ as $n, m \rightarrow \infty$. Let us introduce some more notations that will come in handy in proving this convergence. Given that $\theta \sim \mathbb{P}^N$, we define
\begin{itemize}
    \item $\nu((-\infty, y];t) := P(D(x; \theta) \leq y \mid H(x; \theta) = t)$,
    \item $\hat{\nu}^m((-\infty, y]; t) := P(\hat{D}^m (x; \theta) \leq y \mid \hat{H}^m (x; \theta) = t)$,
\end{itemize}
as the probability measures corresponding to the given conditional distributions. These measures will be useful for characterizing $g(t; x)$ and $\hat{g}^m(t; x)$ respectively where $\hat{g}^m(t; x) = \E_{\theta, \mathbb{P}_{\theta}}[\partial_x \hat{H}^m(x; \theta) \mid \hat{H}^m(x; \theta) = t]$. In what follows, we let $\mathcal{B}_{\eta}(y)$ denote a ball centered at $y$ with radius $\eta$.

\begin{assumption} \label{measure-assumption}
    Assume that there exists a family of measures $G_m(\cdot)$ and a number $\eta > 0$ such that for all $t \in \mathcal{B}_{\eta}(v_{\alpha}(x))$ and for all $\Delta y \subset (-\infty, \infty)$,
    \begin{equation*}
        |\nu(\Delta y, t) - \hat{\nu}^m (\Delta y, t)| \leq G_m(\Delta y) \text{ and } \int_{\mathbb{R}} |y| G_m(d y) \rightarrow 0 \text{ as } m \rightarrow \infty.
    \end{equation*}
\end{assumption}

\begin{assumption} \label{supremum-assumption}
$\sup_{\theta} \E_{\mathbb{P}_{\theta}}[d(x; \xi)^2] < \infty$.
\end{assumption}

Assumptions \ref{measure-assumption} and \ref{supremum-assumption} impose technical conditions to ensure that the estimation errors for both the function value and its gradients are well behaved. Assumption \ref{measure-assumption} is seemingly abstract and deserves further explanation. It essentially requires that the distribution of the gradient estimate conditioned on the function value converges to its true counterpart. One would notice that conditioned on the value of the function estimate, the gradient estimator is no longer unbiased and $\E_{\theta, \mathbb{P}_{\theta}}[\hat{D}^m(x; \theta) - D(x;\theta) \mid \hat{H}^m(x; \theta)=t] \neq 0$ in general, as the observations ($\hat{H}^m$ and $\hat{D}^m$) rely on the same set of $\xi$'s. Intuition suggests that as $m \rightarrow \infty$ and the estimation error $\hat{H}^m(x; \theta) - H(x; \theta) \rightarrow 0$, the corresponding errors in gradient estimation should also cancel out. Assumption \ref{measure-assumption} is a technical condition that we impose to mitigate the difficulties arising from conditioning on measure zero sets in proving this behavior. In fact, if the condition $H(x; \theta) = t$ (and its noisy counterpart) is relaxed from a point to a neighborhood, i.e. $H(x; \theta) \in \mathcal{B}_{\eta}(t)$, it can be shown that Assumption \ref{measure-assumption} follows from Assumptions \ref{zhu-assumption} \& \ref{supremum-assumption}. 
\edit{We provide a detailed discussion on the assumptions in the online supplement, where it is also shown that Assumption \ref{measure-assumption} is satisfied in a general class of problems.}

Now that we have established the necessary regularity conditions, we have the following proposition on the asymptotic bias of $\varphi^{n,m}_{\alpha}(x)$. In the following $a_n = \mathcal{O}(b_n)$ means $\limsup_{n \rightarrow \infty} |a_n / b_n| < \infty$, $a_n = o(b_n)$ means $\lim_{n \rightarrow \infty} a_n/b_n = 0$, and $a_n = \Theta(b_n)$ means $a_n = \mathcal{O}(b_n)$ and $b_n = \mathcal{O}(a_n)$.

\begin{proposition} \label{var-unbiased-proposition}
    Suppose that Assumptions \ref{zhu-assumption} - \ref{supremum-assumption} hold. Then $\E_{\theta, \mathbb{P}_{\theta}}[\varphi^{n,m}_{\alpha}(x)] - v'_{\alpha}(x) \rightarrow 0$ as $n,m \rightarrow \infty$.
    
    Moreover, if in addition the integral in Assumption \ref{measure-assumption} is $\mathcal{O}(m^{-1/2})$, $g(t; x)$ is differentiable w.r.t. $t$ at $t = v_{\alpha}(x)$, and the budget sequence is such that $n = \Theta(m)$, then the bias is $\E_{\theta, \mathbb{P}_{\theta}}[\varphi^{n,m}_{\alpha}(x)] - v'_{\alpha}(x) = \mathcal{O}(n^{-1/2})$.
\end{proposition}

Even though $\varphi^{n,m}$ is asymptotically unbiased, it is not consistent in general when $\theta$ is multidimensional, particularly when the set $\{ \theta: H(x; \theta) = v_{\alpha}(x) \}$ is not a singleton. See \cite{Hong2009VaR} for a discussion on consistency of $I^n$, and \cite{Jiang2015Quantile} for an additional assumption under which $I^n$ is consistent along with some examples. The same argument carries on to our case. A common approach is to use batching to address this difficulty. We have the following theorem which provides the consistency of the batch-mean estimator $\bar{\varphi}^{n,m,k}_{\alpha}(x) := \frac{1}{k} \sum_{i=1}^k \varphi^{n,m}_{\alpha, i}(x)$, where $\varphi^{n,m}_{\alpha, i}(x)$ are i.i.d. copies and k is the number of batches.

\begin{theorem} \label{var-consistency-theorem}
    Suppose that Assumptions \ref{zhu-assumption} - \ref{supremum-assumption} hold, then 
    \begin{equation*}
        \bar{\varphi}^{n,m,k}_{\alpha}(x) \xrightarrow{P} v'_{\alpha}(x) \text{ as } n, m, k \rightarrow \infty,
    \end{equation*}
    where $\xrightarrow{P}$ denotes convergence in probability.
\end{theorem}

In addition to the asymptotic unbiasedness and consistency, we have the following result that characterizes the asymptotic distribution of $\bar{\varphi}^{n,m,k}_{\alpha}(x)$. The proof is a direct application of Lyapunov's Central Limit Theorem combined with Proposition \ref{var-unbiased-proposition}. It is identical to the proof of Theorem 5 of \cite{Hong2009VaR}, and is omitted here.

\begin{theorem} \label{var-normality-theorem}
    Suppose that the (stronger) assumptions of Proposition \ref{var-unbiased-proposition} hold, $k = o(n)$, and $\sup_{n,m}\E_{\theta, \mathbb{P}_{\theta}}[|\varphi^{n,m}|^{2 + \gamma}] < \infty$ for some $\gamma > 0$. Then,
    \begin{equation*}
        \sqrt{k}(\bar{\varphi}^{n,m,k} - v'_{\alpha}) \Rightarrow \mathcal{N}(0, \sigma_{\infty}) \ \text{as} \ n,m,k \rightarrow \infty,
    \end{equation*}
    where $\sigma^2_{\infty} = \lim_{n,m \rightarrow \infty} Var(\varphi^{n,m})$ is the asymptotic variance of $\varphi^{n,m}$.
\end{theorem}

\subsubsection{Conditional Value-at-Risk Case} 
Conditional Value-at-Risk, defined as $\text{CVaR}_{\alpha} = \E_{\theta} [H(x; \theta) \mid H(x; \theta) \geq v_{\alpha}(x)]$, is the expectation of large losses. We are interested in estimating the gradient $d \text{CVaR}_{\alpha}(H(x; \theta)) / d x$ using samples of $h(x, \xi(\theta))$ (and $d(x, \xi(\theta))$). We use $c_{\alpha}(x)$ and $c'_{\alpha}(x)$ as shorthand notations for CVaR$_{\alpha}(H(x; \theta))$ and $d \text{CVaR}_{\alpha}(H(x; \theta)) / dx$ respectively.

Under a mild set of assumptions, \cite{Hong2009CVaR} show that CVaR gradients can be written in the form of a conditional expectation as 
\begin{equation}\label{cvar-expectation}
    c'_{\alpha}(x) = \E_{\theta}[D(x; \theta) \mid H(x; \theta) \geq v_{\alpha}(x)].
\end{equation}
We propose the following estimator of CVaR gradients that mimics a Monte-Carlo estimator of (\ref{cvar-expectation}) with the available information.
\begin{equation}
    \psi^{n,m}_{\alpha}(x) := \frac{1}{n(1-\alpha)} \sum_{i=1}^n \hat{D}^m (x; \theta_i) \mathds{1}_{(\hat{H}^m (x; \theta_i) \geq \hat{v}^{n,m}_{\alpha}(x) )}.
\end{equation}
In the remainder of this subsection, we show that $\psi^{n,m}_{\alpha}(x)$ is a strongly consistent and asymptotically unbiased estimator of $c'_{\alpha}(x)$.

The analysis of $\psi^{n,m}_{\alpha}(x)$ relies on a weaker set of assumptions than that of $\varphi^{n,m}_{\alpha}(x)$. Assumption \ref{measure-assumption} is no longer needed, and Assumption \ref{hong-var-assumptions} is replaced with the following weaker assumption due to \cite{Hong2009CVaR}. Assumption \ref{hong-assumption}, along with Assumption \ref{interchange-assumption}, is needed to ensure validity of (\ref{cvar-expectation}).

\begin{assumption} \label{hong-assumption}
\citep{Hong2009CVaR}
\begin{enumerate}
    \item The VaR function $v_{\alpha}(x)$ is differentiable for any $x \in \mathcal{X}$.
    
    \item For any $x \in \mathcal{X}$, $P[H(x; \theta) = v_{\alpha}(x)] = 0$.
\end{enumerate}
\end{assumption}

\edit{We note that Assumption \ref{hong-assumption} and is implied by \ref{hong-var-assumptions} and the differentiability of $h(x, \xi(\theta))$. It is presented separately here, as it replaces Assumption \ref{hong-var-assumptions} with a weaker set of conditions.}
The following proposition is needed in proving the consistency of $\psi^{n,m}_{\alpha}(x)$. 

\begin{proposition} \label{indicator-proposition}
Suppose Assumption \ref{zhu-assumption} holds and $P(H(x; \theta) = v_{\alpha}(x)) = 0$. Then
$$ \frac{1}{n} \sum_{i=1}^n | \mathds{1}_{(\hat{H}^m (x; \theta_i) \geq \hat{v}^{n,m}_{\alpha}(x))} - \mathds{1}_{(H(x; \theta_i) \geq v_{\alpha}(x))} | \rightarrow 0 \ w.p.1 \text{ as } n,m \rightarrow \infty .$$
\end{proposition}

Note that (\ref{cvar-expectation}) can be rewritten as $c'_{\alpha}(x) = \frac{1}{1-\alpha} \E_{\theta, \mathbb{P}_{\theta}}[d(x, \xi) \mathds{1}_{(H(x; \theta) \geq v_{\alpha}(x))}]$, which admits $\frac{1}{n (1-\alpha)} \sum_{i=1}^n \hat{D}^m(x; \theta) \mathds{1}_{(H(x; \theta) \geq v_{\alpha}(x))}$ as a Monte-Carlo estimator. Proposition \ref{indicator-proposition} shows that the bias introduced by replacing $\mathds{1}_{(H(x; \theta) \geq v_{\alpha}(x))}$ with its noisy version, $\mathds{1}_{(\hat{H}^m (x; \theta_i) \geq \hat{v}^{n,m}_{\alpha}(x) )}$, disappears in the limit. The following proposition extends this result to show that $\psi^{n,m}_{\alpha}(x)$ is asymptotically unbiased and the bias is of the order $\mathcal{O}(n^{-1/2})$.

\begin{proposition} \label{cvar-bias-proposition}
Under Assumptions \ref{zhu-assumption}, \ref{interchange-assumption}, \ref{supremum-assumption}, \ref{hong-assumption}, the bias $\E_{\theta, \mathbb{P}_{\theta}}[\psi^{n,m}_{\alpha}(x)] - c'_{\alpha}(x) \rightarrow 0$ as $n,m \rightarrow \infty$.
Moreover, if in addition $n = \Theta(m)$, then the bias is of the order $\mathcal{O}(n^{-1/2})$.
\end{proposition}

We conclude this subsection with the following theorem that provides strong consistency of $\psi^{n,m}_{\alpha}(x)$.

\begin{theorem} \label{consistency-theorem}
Under Assumptions \ref{zhu-assumption}, \ref{interchange-assumption}, \ref{supremum-assumption}, \ref{hong-assumption}, $\psi^{n,m}_{\alpha}(x)$ is a strongly consistent estimator of $c'_{\alpha}$. 
\end{theorem}

\begin{remark}
Even though the results here are derived using the IPA estimator for the inner expectation, one should notice that our proofs only require a consistent estimator of the inner expectation. Therefore, where IPA is not applicable or it is not preferred for any other reason, one could replace $\hat{D}^m$ with any consistent estimator such as the generalized likelihood ratio estimator of \cite{Peng2018GLR}, support independent unified likelihood ratio and infinitesimal perturbation analysis estimator of \cite{Wang2012SLRIPA} etc. as long as the corresponding regularity conditions hold.
\end{remark}

\subsection{Convergence Analysis of the Algorithms}
\edit{
In this subsection, we start with a brief discussion on implementation and computational cost of the SA algorithm, and show that the use of $\varphi^{n,m}_{\alpha}(x)$ and $\psi^{n,m}(x)$ results in consistent algorithms for solving the corresponding BRO problems. 
}

\edit{
The SA algorithm is briefly summarized in Algorithm \ref{alg:step-by-step}. When a closed form of the posterior distribution is not available, one may use numerical methods, such as Markov Chain Monte Carlo (MCMC) \citep{Lange2010}, variational Bayes \citep{Fox2012VariationalBayes} etc., to approximate the posterior distribution. 
We emphasize that it is only necessary to draw an empirical approximation to the posterior before the optimization starts (see Section \ref{two-sided-section} for more details). This avoids a repeated use of e.g. MCMC, which is not necessary since the posterior distribution does not change, and facilitates cost effective sampling of $\theta \sim \mathbb{P}^N$ from the generated empirical distribution. Thus, the computational cost of Algorithm \ref{alg:step-by-step} is dominated by the simulations of $h(\cdot, \cdot)$, which has a total computational cost of $\mathcal{O}(n_T m_T T^2)$.
}

\edit{
\begin{algorithm}[tb]
  \caption{Optimization of BRO via SA}
  \label{alg:step-by-step}
\begin{algorithmic}
  \STATE {\bfseries Input:} \edit{$\{n_t\}_{t \geq 0}, \{m_t\}_{t \geq 0}, \{\epsilon_t\}_{t \geq 0} x_0, T$, the input data, prior distribution and the choice of risk measure.}
  \STATE \edit{Calculate the posterior distribution $\mathbb{P}^N$. If a closed form or $\mathbb{P}^N$ is not available, use a computational method (e.g. MCMC or variational Bayes) to draw an empirical approximation.}
  \FOR{$t=0$ {\bfseries to} $T-1$}
    \STATE \edit{ Draw $\theta_1, \ldots, \theta_{n_t} \stackrel{iid}{\sim} \mathbb{P}^N$ and $\xi_1(\theta_i), \ldots, \xi_{m_t}(\theta_i) \stackrel{iid}{\sim} \mathbb{P}_{\theta_i}$.}
    \STATE \edit{Simulate $h(x, \xi_j(\theta_i))$, and calculate the estimator $\varphi_{\alpha}^{n_t, m_t}$ if $\rho$ is chosen as VaR, or $\psi_{\alpha}^{n_t,m_t}$ if $\rho$ is chosen as CVaR.}
    \STATE \edit{ Set $x_{t+1} = \Pi_{\mathcal{X}}[x_t + \epsilon_t Y_t]$, where $Y_t$ is either of $\varphi_{\alpha}^{n_t, m_t}$ if $\rho$ is chosen as VaR, or $\psi_{\alpha}^{n_t,m_t}$ if $\rho$ is chosen as CVaR.}
  \ENDFOR
  \STATE {\bfseries Return:} \edit{$x_T$ as the decision.}
\end{algorithmic}
\end{algorithm}
}

\edit{
The remainder of the subsection is dedicated to proving the convergence of the algorithms.
The main result is summarized in the following theorem.
}

\begin{theorem} \label{rm-convergence-theorem}
Suppose that Assumptions \ref{zhu-assumption} - \ref{supremum-assumption} hold, and $v_{\alpha}(\cdot)$ is continuously differentiable w.r.t $x$, $\{n_t\}, \{m_t\}$ are monotonically increasing sequences, $\sum_{t=0}^{\infty} \epsilon_t = \infty$, $\sum_{t=0}^{\infty} \epsilon^2_t < \infty$. Then, the SA algorithm (\ref{rm-algorithm}) with $Y_t = -\bar{\varphi}^{n_t, m_t, k_t}_{\alpha}(x_t)$ converges w.p.1 to a unique solution set of the ODE
\begin{equation}
    \dot{x} = -v'_{\alpha}(x).
\end{equation}

Similarly, under Assumptions \ref{zhu-assumption}, \ref{interchange-assumption}, \ref{supremum-assumption}, \ref{hong-assumption}, and assuming $c_{\alpha}(\cdot)$ is continuously differentiable w.r.t $x$, $\{n_t\}, \{m_t\}$ are monotonically increasing sequences, $\sum_{t=0}^{\infty} \epsilon_t = \infty$, $\sum_{t=0}^{\infty} \epsilon^2_t < \infty$; the SA algorithm (\ref{rm-algorithm}) with $Y_t = -\psi^{n_t, m_t}_{\alpha}(x_t)$ converges w.p.1 to a unique solution set of the ODE
\begin{equation}
    \dot{x} = -c'_{\alpha}(x).
\end{equation}
\end{theorem}

\edit{
Theorem \ref{rm-convergence-theorem} is a direct application of Theorem 2.1 of \cite{Kushner&Yin2003SA}.
Following their analysis, $Y_t$ is deconstructed as $Y_t = g(x_t) + \delta M_t + \beta_t$, where $g(x_t)$ is the negative gradient at $x_t$, $\delta M_t$ is martingale difference error term, and $\beta_t$ is the bias term. \cite{Kushner&Yin2003SA} imposes the following set of assumptions to ensure the convergence of the SA algorithm.
}

\newpage

\begin{assumption} Chapter 5.2, \cite{Kushner&Yin2003SA} \label{kushner_assumptions}
\begin{enumerate}
    \item $\sup_t \E[Y_t^2] < \infty$;
    \item There is a measurable function $g(\cdot)$ of $x$ and random variables $\beta_t$ such that 
    $$\E[Y_t \mid x_0, Y_i, i<t] = g(x_t) + \beta_t ;$$
    \item $g(\cdot)$ is continuous;
    \item $\sum_{t=0}^{\infty} \epsilon_t = \infty$, $\sum_{t=0}^{\infty} \epsilon_t^2 < \infty$;
    \item $\sum_{t=0}^{\infty} \epsilon_t |\beta_t| < \infty$ w.p.1.
\end{enumerate}
\end{assumption}

Let $Y_t = -\bar{\varphi}^{n_t, m_t, k_t}_{\alpha}(x_t)$ for the VaR estimator with $g(x_t) = -v'_{\alpha}(x_t)$, $\delta M_t = \E[\bar{\varphi}^{n_t,m_t,k_t}_{\alpha}(x_t)] - \bar{\varphi}^{n_t,m_t,k_t}_{\alpha}(x_t)$, and $\beta_t = v'_{\alpha}(x_t) - \E[\bar{\varphi}^{n_t,m_t,k_t}_{\alpha}(x_t)]$. 
For the CVaR, apply the same decomposition with $\psi^{n_t, m_t}_{\alpha}(x_t)$ replacing $\bar{\varphi}^{n_t,m_t,k_t}_{\alpha}(x_t)$. To see that Assumption \ref{kushner_assumptions} is satisfied, note the following. Assumption \ref{kushner_assumptions}.1 immediately follows from Assumption \ref{supremum-assumption}. Assumption \ref{kushner_assumptions}.2 is satisfied by the given deconstruction. For Assumption \ref{kushner_assumptions}.3, we assume that $v_{\alpha}(\cdot)$ and $c_{\alpha}(\cdot)$ are continuously differentiable. Assumption \ref{kushner_assumptions}.4 is a common requirement for the step size sequences and is imposed here. As shown by Theorem 2.3 of \cite{Kushner&Yin2003SA} and Theorem 2 of \cite{Kushner2010SA-Survey},  Assumption \ref{kushner_assumptions}.5 can be replaced with $\beta_t \rightarrow 0$ w.p.1 which is given by Propositions \ref{var-unbiased-proposition} and \ref{cvar-bias-proposition} for VaR and CVaR cases respectively. Therefore, Assumption \ref{kushner_assumptions} is satisfied and Theorem \ref{rm-convergence-theorem} follows immediately from Theorem 2.1 of \cite{Kushner&Yin2003SA}.

\begin{remark} \label{batching-remark}
    A careful look at the Assumption \ref{kushner_assumptions} reveals that in order for the convergence in Theorem \ref{consistency-theorem} to hold, we do not need a consistent estimator. Therefore, when the optimization is of concern, one can opt to use the non-batching VaR estimator $\varphi^{n,m}_{\alpha}(x)$ without sacrificing the convergence of the algorithm. 
\end{remark}

\section{Numerical Examples} \label{numerical-example-section}
\edit{
In this section, we present an empirical study of the proposed algorithm.
We start with a simple quadratic example, where we compare the numerical efficiency of two gradient-based algorithms that use the estimators developed in this paper with two gradient-free approaches from the existing literature. 
We follow that with a more realistic example of a two-sided market model, where we demonstrate the convergence of the SA algorithm on several BRO objectives. The section is concluded with a discussion on the objective choice, where the robustness of various objective choices are demonstrated. 
}

\edit{
\subsection{A simple quadratic example}
In this section, we study a simple quadratic example, where we compare the numerical efficiency of optimization algorithms using the gradient estimators developed in this paper with the gradient-free methods from the literature that can be used to solve the BRO problem. For gradient-based methods, we consider a quasi-newton method, the LBFGS algorithm \citep{Zhu1997LBFGS}, which only requires access to the gradients of the function, as well as the SA algorithm described above. For gradient-free alternatives, we consider the Nelder-Mead simplex method \citep{NelderMead1965}, and the Expected Improvement algorithm \citep{Jones1998EI}, both of which are known for their superior empirical performance.

The example in consideration is modified from \cite{Hong2009VaR}, and is given by $H(x; \theta) = x \theta_1 + x^2 \theta_2$ with the simulation oracle $h(x, \xi(\theta)) = x \theta_1 + x^2 \theta_2 + x \xi(\theta)$, where $\xi(\theta) \sim \mathcal{N}(0, \frac{\theta_1^2}{100})$. It follows that 
$D(x; \theta) = \theta_1 + 2x \theta_2$, and
$d(x, \xi(\theta)) = \theta_1 + 2x \theta_2 + \xi(\theta)$.

In the online supplement, we discuss the details of the experiment, verify the assumptions, and obtain the analytical solution to the BRO optimization problem. 
Table \ref{simple-results} presents the average optimality gap obtained from $50$ replications using the BRO CVaR objective with risk level $\alpha = 0.75$. The algorithms use the same simulation budget, where the BRO objective and its gradient is estimated using $n_t = n = 100$ and $m_t = n_t/5$. Since the benchmark algorithms are developed for deterministic optimization, we consider both stochastic evaluations of the objective and the Sample Average Approximation (SAA, \citealp{Kim2015SAA}) counterpart, which converts it into an approximate deterministic optimization problem by fixing the random variables $\theta$ and $\xi$. 

\begin{table}[ht]
\centering
\caption{The optimality gap in the simple quadratic example. The reported values are on the scale of $10^{-2}$.}
{\begin{tabular}{cl|rrrr}
\# of evaluations & SAA?   & SA     & LBFGS     & Nelder-Mead     & EI      \\ 
\hline
\multirow{2}{*}{10}  & No      & \colorbox{green}{$1.131$} & $6.114$ & $244.515$ & $6.055$  \\
                      & Yes  & $2.054$  & $0.958$  & $17.146$  & $13.312$   \\ 
\multirow{2}{*}{20}  & No & \colorbox{green}{$0.138$}  & $6.274$  & $156.030$  & $6.164$   \\ 
                      & Yes & $1.057$  & $0.958$  & $0.955$  & $1.607$   \\ 
\multirow{2}{*}{50}  & No & \colorbox{green}{$0.036$}  & $6.274$  & $156.402$  & $4.149$   \\ 
                      & Yes & $0.959$  & $0.958$  & $0.958$  & $1.003$   \\ 
\multirow{2}{*}{100}  & No & \colorbox{green}{$0.015$}  & $6.274$  & $156.401$  & $4.257$   \\ 
                      & Yes & $0.958$  & $0.958$  & $0.958$  & $0.957$   \\ 
\hline
\end{tabular}}
\label{simple-results}
\end{table}

The results show a clear advantage of using the SA algorithm over all the benchmarks considered. Using the stochastic gradient estimators, our proposed Algorithm \ref{alg:step-by-step} (results highlighted in Table \ref{simple-results}) achieves almost 2 orders of magnitude better performance than the closest competitor. We observe that the benchmark algorithms (LBFGS, Nelder-Mead and EI) have difficulty in solving the optimization problem using the stochastic estimators, while the results improve a little when solving the SAA counterpart. This shows that the methods developed in this paper provide a clear improvement over the existing alternatives for optimizing the BRO problem.

}

\subsection{A Two-Sided Market Model} \label{two-sided-section}
In a two-sided market model, the customers and providers arrive to the system according to two independent arrival processes. Upon arrival, a customer is served immediately if there is an available provider, otherwise the customer queues up to be served by future provider arrivals. Similarly, arriving providers leave the queue immediately if there is a customer waiting, otherwise they wait for the future customer arrivals. With some slight variation, such models can be used to mimic system dynamics of various real life scenarios such as sharing or gig economies. In this example, it is assumed that a provider can only serve one customer, and the system operates without abandonment. Customer arrivals are assumed to follow a Poisson process with rate $\lambda(p)$, and provider arrivals follow a Poisson process with rate $\mu(p)$, where $p$ denotes the price set by the platform, with $\lambda(p)$ ($\mu(p)$) decreasing (increasing) in $p$. The rate functions are given by $\lambda(p) = K^C \frac{2 exp(- \theta^C p)}{1 + exp(- \theta^C p)}, \mu(p) = K^P \frac{1 - exp(- \theta^P p)}{1 + exp(- \theta^P p)}$ where $K^C, K^P$ are the (known) potential numbers of customers and providers, and $\theta = (\theta^C, \theta^P)$ are the (unknown) sensitivities of customers and providers respectively. These rate functions result in $\lambda(0) = K^C, \lim_{p \rightarrow \infty} \lambda(p) = 0, \mu(0) = 0$ and $\lim_{p \rightarrow \infty} \mu(p) = K^P$, which agrees with the intuition that no providers (customers) should be willing to participate when the price is $0$ ($\infty$) and vice versa.

In our setting, the platform aims to minimize customer wait time to improve customer satisfaction. However, one could easily see that a naive objective of minimizing the expected wait time would drive price to infinity, leading to excess number of providers and no customers, thus no service or revenue. To avoid this pitfall, the objective is modified to be a weighted combination of customer waiting time and expected revenue. We estimate the customer waiting time by $\frac{1}{M} \sum_{i=1}^M W_i$, the average waiting time of first $M$ customers, where $W_i$ denotes the waiting time of the $i^{th}$ customer, and expected revenue is estimated as $p \lambda(p)$. The resulting objective takes the form
\begin{equation*}
    \min_{p} H(p; \theta) = \E\left[ \frac{1}{M} \sum_{i=1}^M W_i - a p \lambda(p) \right],
\end{equation*}
where $a$ is a predetermined weight. If the rate functions $\lambda(p)$ and $\mu(p)$ (i.e. $\theta$) are known, one can use sampling to estimate and optimize the objective. Since $\theta$ is unknown and is estimated from a finite set of real world data, the objective is replaced by $\min_{p} \rho_{\theta} \{H(p; \theta)\}$ to account for input uncertainty. 

In order to estimate the gradient of the objective, we need to sample from $h(p, \xi(\theta)) = \frac{1}{M} \sum_{i=1}^M W_i - ap\lambda(p)$ and its gradient $dh(p, \xi(\theta)) / dp = \frac{1}{M} \sum_{i=1}^M dW_i/dp - ad(p\lambda(p))/dp$. Note that $W_i = \max \{0, A^P_i - A^C_i \}$, where $A^P_i$ and $A^C_i$ denote the arrival time of $i^{th}$ provider and $i^{th}$ customer respectively, and $dW_i/dp = \mathds{1}_{\{W_i > 0\}} (dA^P_i / dp - dA^C_i / dp)$ where $\frac{dA^P_i}{dp} = \frac{dA^P_i}{d \mu(p)} \frac{d \mu(p)}{dp}$ and $\frac{dA^C_i}{dp} = \frac{dA^C_i}{d \lambda(p)} \frac{d \lambda(p)}{dp}$. The gradient $d(p\lambda(p))/dp$ can be calculated as $d(p\lambda(p))/dp = \lambda(p) + p \frac{d \lambda(p)}{dp}$.

Before we can run the experiments, we need to estimate the objective function $\rho_{\theta} \{H(p; \theta)\}$ (and its gradient) which requires sampling from $\theta \sim \mathbb{P}^N$. One should notice that regardless of the choice of prior, $\mathbb{P}^N$ does not admit a simple closed form solution. However, we can use an MCMC method, the Metropolis-Hastings algorithm (see \citealp{Lange2010}), to sample from $\mathbb{P}^N$. Suppose that the true parameters are $K^C = 40, K^P = 20, \theta^C = 0.1, \theta^P = 0.05$, and we are given a dataset $\phi^N = \{\phi^C, \phi^P \}$ of size $N = 10$ (each) of inter-arrival times drawn at $p = 10$. 

Since the likelihood functions of $\theta^C$ and $\theta^P$ are separable, we estimate the posteriors using two independent MCMC runs. For the MCMC, we use a Gaussian proposal distribution, and the candidate is generated as $\theta_{\text{candidate}} = \theta_{\text{current}} + \mathcal{N}(0, \sigma^2)$ with $\sigma = 2.5 \times 10^{-2}$. We use $\text{Uniform} (0.01, 0.5)$ as an uninformative prior. With the given choice of proposal and prior distributions, the acceptance probability simplifies to 
\begin{equation*}
    P(\text{accept}) = \min \{1, \mathds{1}_{\text{candidate} \in (0.01, 0.5)} p(\phi^N \mid \theta_{\text{candidate}}) / p(\phi^N \mid \theta_{\text{current}}) \},
\end{equation*}
where the likelihood of $\phi^N = \{\xi_1, \ldots, \xi_N\}$ is calculated as $p(\phi^N \mid \theta) = \prod_{i=1}^N f(\xi_i \mid \theta)$ with $f(\cdot \mid \theta)$ as the probability density given the parameter $\theta$. We use a (post burn-in) run length of $10^6$ iterations with the starting point of $\theta_0 = 0.075$ and a burn-in period of $10^5$ iterations.
An empirical analysis of the output using the Wasserstein distance, \edit{described in the online supplement,} suggests that the samples are drawn from a stationary distribution.
Let $\Tilde{\Theta}$ denote the list of $10^6$ samples generated from the MCMC run. $\theta \sim \mathbb{P}^N$ is sampled as follows: we generate a random variable $i \sim \text{discrete-uniform}[1, 10^6]$ as the index and set $\theta = \Tilde{\Theta}[i]$. Since the MCMC converges to the posterior distribution and the $\tilde{\Theta}$ are samples from the approximate steady state distribution of the MCMC, the $\theta$ generated this way are approximately distributed as $\mathbb{P}^N$. The resulting samples from MCMC have a sample average of $5.2 \times 10^{-2}$ and sample standard deviation of $3.2 \times 10^{-2}$ for $\theta^C$, and a sample average of $6.8 \times 10^{-2}$ and sample standard deviation of $2.3 \times 10^{-2}$ for $\theta^P$. The corresponding maximum likelihood estimators (MLE) are given by $\hat{\theta}^C_{\text{MLE}} = 6.06 \times 10^{-3}$ and $\hat{\theta}^P_{\text{MLE}} = 5.9 \times 10^{-2}$, which suggests that the input data for the customers is not representative of the true distribution.

The problem parameters are set as $M=100, a=1/25$. Before going into the optimization, we need to pick the budget and step size sequences. 
\edit{The choice of step size is problem specific, as both too large and too small step size sequences harm the convergence of the algorithm. We recommend using simple pilot experiments to guide the selection.}
For the budget sequences, Propositions \ref{var-unbiased-proposition} \& \ref{cvar-bias-proposition} suggest that $n$ and $m$ should be of the same order, however the relative magnitudes are again problem dependent. We recommend checking the relative magnitudes of stochastic and input uncertainties, estimated by the standard deviation of $\hat{H}^m(p, \theta)$ for a fixed $\theta$ and w.r.t. $\theta$ with $\theta \sim \mathbb{P}^N$ respectively, and modifying $n$ and $m$ until a balance is achieved. As a result of pilot experiments, we pick $m = n/10$ for this example. In light of Remark \ref{batching-remark}, we use the non-batching estimator $\varphi^{n,m}$ for VaR. 

For the algorithm runs, the step size sequence is chosen as $\epsilon_t = \frac{20}{(100+t)^{0.8}}$ and the budget sequence is $n_t = 100 + 0.5 t, m_t = \lfloor n_t / 10 \rfloor$. For each choice of $\rho$ below, we run $50$ replications of the algorithms, each for $1000$ iterations with $p_0 = 5$. The results are reported in Table \ref{offline-results}. For each $\rho$, we report the average solution obtained from $50$ replications (as $p$), the estimate of the solution standard deviations (as std($p$)), the approximate optimal solution to the corresponding BRO problem (as $p^*_{\rho}$), and the performance of the obtained solution under the true distribution (as $H^c(p)$). $H^c(p)$ was evaluated using Monte Carlo simulation with $10^5$ samples. The reported $p^*_{\rho}$ is computed via brute-force Monte Carlo simulation with simulation intervals of $0.1$ using common random numbers (CRN) and a budget of $n = 10^4, m = 10^3$. One would notice that for certain choices of $\rho$ there is a discrepancy of about $0.1 - 0.2$ between the algorithm solutions and the estimated optimal solutions. We note that the difference between the solution performances was below $10^{-2}$ in each case, and the difference can be attributed to the estimator bias. For comparison, the true optimal solution and its performance is estimated as $p^*_c = 20.47$ and $H^c(p^*_c) = -7.160$, using Monte Carlo simulation with simulation intervals of $0.01$ using CRN and $4 \times 10^5$ samples. The MLE solution is estimated in a similar manner to be $p^*_{\text{MLE}} = 211$ with $H^c(p^*_{\text{MLE}}) = -4.63 \times 10^{-7}$, which points to the value of robustness in this particular example.

\begin{table}[ht]
\centering
\caption{Algorithm solutions, solution standard deviations, the approximate optimal solutions, and the solution performance under the true problem.}
{\begin{tabular}{cl|rrrrr}
\multicolumn{2}{c}{$\alpha$}   & 0.5     & 0.6     & 0.7     & 0.8     & 0.9      \\ 
\hline
\multirow{4}{*}{VaR}  & $p$      & $31.026$ & $27.423$ & $24.471$ & $22.118$ & $19.969$  \\
                      & std($p$) $(\times 10^{-1})$ & $1.24$  & $1.20$  & $1.31$  & $1.73$  & $2.01$   \\ 
                      & $p^*_{\rho}$ & $30.9$  & $27.3$  & $24.3$  & $21.9$  & $19.9$   \\ 
                      & $H^c(p)$ & $-4.269$  & $-5.309$  & $-6.230$  & $-6.924$  & $-4.561$   \\ 
\hline
\multirow{4}{*}{CVaR} & $p$      & $23.365$ & $22.049$ & $20.898$ & $19.868$ & $18.904$  \\
                      & std($p$) $(\times 10^{-2})$ & $1.97$  & $1.65$  & $1.76$  & $2.37$  & $2.73$   \\
                      & $p^*_{\rho}$ & $23.3$  & $22.0$  & $20.9$  & $19.8$  & $18.8$   \\ 
                      & $H^c(p)$ & $-6.573$  & $-6.941$  & $-7.142$  & $-7.123$  & $-6.914$   \\ 
\hline
\end{tabular}}
\label{offline-results}
\end{table}

\begin{figure}[ht]
    \centering
    \includegraphics[height = 4in]{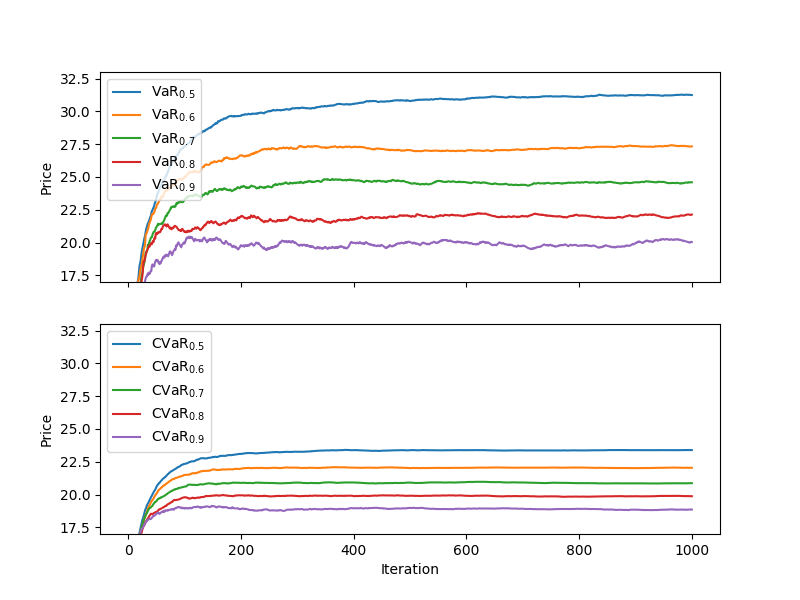}
    \caption{The evolution of the solution (Price) through iterations in a typical algorithm run.}
    \label{fig:offline-plots}
\end{figure}

In Figure \ref{fig:offline-plots}, we plot a typical algorithm run for each $\rho$. It is seen that the algorithm solutions quickly move into a neighborhood of the optimal solution and proceed to refine the solution further in the following iterations. There is a striking difference between the variability of the solution paths corresponding to VaR and CVaR objectives, which can be attributed to our decision to forgo batching in favor of computational efficiency. Without batching, the value of the VaR estimator is calculated using a single realization of $\theta$, whereas, the value of the CVaR estimator is calculated by averaging over a number of $\theta$'s. 

The results reported in Table \ref{offline-results} demonstrate the convergence of the algorithms to a the optimal solution as given in Theorem \ref{rm-convergence-theorem}. Moreover, it is seen that different choices of $\rho$ correspond to a wide spread of solutions, which in turn has a significant effect on the resulting objective values. The true performance (i.e., under the true input parameter) of the solutions will be shown in the next section, while we discuss the choice of $\rho$.

\subsection{Discussion on objective choice} \label{objective-choice-section}
So far, our work focuses on solving the BRO problem given a risk function. However, the choice of the risk function (or the objective) is not a trivial task by itself. In this subsection, we will empirically compare several objectives and try to highlight the effect each has on the resulting decision. We draw $50$ independent input data sets of size $N = 10$ each. For each set of input data, we estimate the posterior distribution and optimize the corresponding objective functions using the same parameters as the original problem. For VaR and CVaR objectives, we use risk levels $\alpha \in \{0.5, 0.7, 0.9\}$. In addition, we compare with the Expectation, Mean-Variance (with variance weight of $0.1$) and MLE objectives.

\begin{figure}[ht]
    \centering
    \includegraphics[width=5.5in]{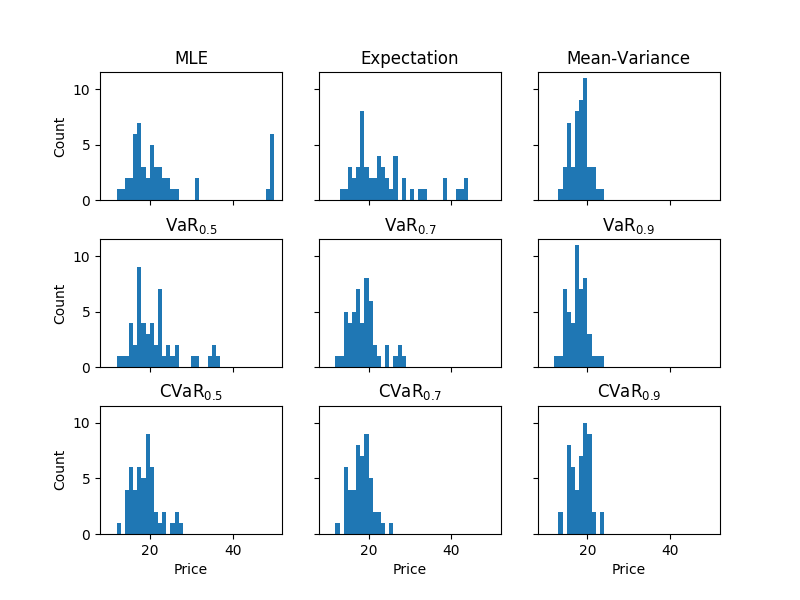}
    \caption{\edit{Optimal solutions obtained from various choices of $\rho$.}}
    \label{fig:rho_x_comparison}
\end{figure}

\begin{figure}[ht]
    \centering
    \includegraphics[width=5.5in]{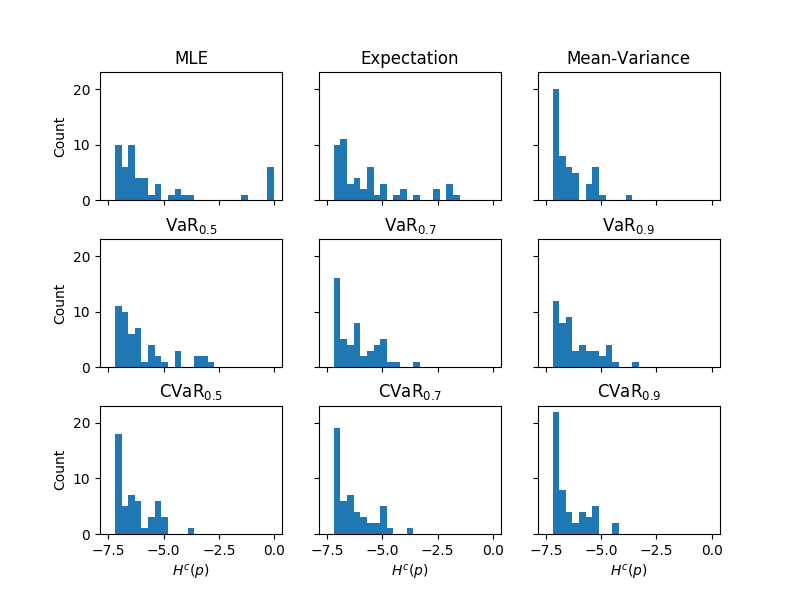}
    \caption{\edit{Performance (under true input parameter) of optimal solutions obtained from various choices of $\rho$.}}
    \label{fig:rho_comparison}
\end{figure}

The performance of the solutions obtained from algorithm runs are estimated via Monte Carlo simulation using CRN and $10^4$ samples. The histograms of the solutions and the solution performances are then plotted in Figures \ref{fig:rho_x_comparison} and \ref{fig:rho_comparison} for each choice of $\rho$. In order to highlight the important areas, we restrict the histograms to a range of $[10, 50]$ in Figure \ref{fig:rho_x_comparison}. This was only an issue for the case of MLE, where the solutions ranged up to $500$. Any solution value that exceeded $50$ is plotted as a point at $50$. 

A quick look at Figures \ref{fig:rho_x_comparison} \& \ref{fig:rho_comparison} reveals the importance of objective choice. Figure \ref{fig:rho_x_comparison} demonstrates the robustness of the BRO objectives, in the sense that the solutions are robust to the particular realization of the input data, and are more concentrated compared to MLE objectives. Although not explicitly shown, observation of just a few outliers in the input data affects the resulting MLE solutions drastically, whereas, the risk averse BRO solutions are much less sensitive to such observations. The choice of a small input size of $N = 10$ further highlights the importance of the objective choice here. As expected, the robustness increases with $\alpha$ and the CVaR objectives are more robust than VaR objectives by definition. In this example, the BRO solutions tend to concentrate around the true optimal solution, which leads to superior overall performance compared to MLE objective, as seen in Figure \ref{fig:rho_comparison}. We would like to emphasize that, although preferred, the superior solution performance is not something that a robust objective aims to provide, and the true aim of a robust objective is to provide a consistent solution performance across a wide range of input data. We refer an interested reader to \cite{Zhou2017book-chapter} for a similar numerical study on an M/M/1 queue problem and a News-vendor problem. We end our discussion by noting that it is possible to combine several objectives studied here and solve them using the tools developed in this paper. For example, if one wishes to balance between the robustness of VaR \& CVaR and the average solution performance, Mean-VaR and Mean-CVaR objectives are obvious choices. In addition to choosing $\alpha$, one can adjust the relative weights of the Mean and VaR/CVaR objectives to balance between robustness and expected solution performance.

\section{Conclusion} \label{conclusion}
In this paper, with the aim of developing efficient methods for solving the Bayesian Risk Optimization framework, we derive stochastic gradient estimators and propose associated stochastic approximation algorithms. Our estimators extend the literature of stochastic gradient estimation to the case of nested risk functions. An example of a two-sided market model is studied to demonstrate the numerical performance of the algorithms, and provide insight into the choice of BRO objectives. Although the exposition of the paper focuses on the BRO framework, the gradient estimators we develop can be used in other settings where nested simulation is used, such as estimating the sensitivities of complex financial portfolios.

\section*{Acknowledgements}
The authors gratefully acknowledge the support by the National Science Foundation under Grant CAREER
CMMI-1453934 and the Air Force Office of Scientific Research under Grant FA9550-19-1-0283.

\bibliographystyle{apacite} 
\bibliography{ref.bib}

\end{document}


\title{Online Supplement to ``Solving Bayesian Risk Optimization via Nested Stochastic Gradient Estimation" by Cakmak, Wu and Zhou}

\date{}

\maketitle

\tableofcontents

\section{Overview}

Recall that our estimators are given by
\begin{equation*}
    \varphi^{n,m}_{\alpha}(x) := \partial_x \hat{H}^m (x; \theta) |_{ \hat{H}^m (x; \theta) = \hat{v}^{n,m}_{\alpha} } = \hat{D}^m (x; \hat{\theta}^m_{(\lceil \alpha n \rceil)})
\end{equation*}
for VaR; and 
\begin{equation*}
    \psi^{n,m}_{\alpha}(x) := \frac{1}{n(1-\alpha)} \sum_{i=1}^n \hat{D}^m (x; \theta_i) \mathds{1}_{(\hat{H}^m (x; \theta_i) \geq \hat{v}^{n,m}_{\alpha}(x) )}
\end{equation*}
for CVaR. 

\edit{
For the readers' convenience, we repeat the full set of assumptions here. We continue with a detailed discussion on the assumptions, verify them on a simple example, and show that Assumption 3.5 is satisfied for a general class of problems. 
We fill in the details of the numerical experiments, and conclude with the proofs of the results presented in the paper.
}

\section{Complete list of assumptions}

\begin{manualassumption}{3.2}  \label{zhu-assumption}
\cite{Zhu2017input-uncertainty}
\begin{enumerate}
    \item For all $x \in \mathcal{X}$, the response $h(x, \xi(\theta))$ has finite conditional second moment, i.e., \linebreak $\tau_{\theta}^2 = \mathbb{E}_{\mathbb{P}_{\theta}}[h(x, \xi)^2] < \infty$ w.p. 1 $(\mathbb{P}^N)$ and $\tau^2 = \mathbb{E}_{\theta, \mathbb{P}_{\theta}}[h(x, \xi)^2] = \int \tau_{\theta}^2 d\mathbb{P}^N < \infty$.
    
    \item The joint density $p_m(h,e)$ of $H(x; \theta)$ and $\bar{\mathcal{E}}^m(x; \theta)$, and its partial gradients $\frac{d}{dh} p_m(h,e)$ and $\frac{d^2}{dh^2}p_m(h,e)$ exist for each $m$, all pairs of $(h,e)$ and for all $x \in \mathcal{X}$.
    
    \item For all $x \in \mathcal{X}$, there exists non-negative functions $g_{0,m}(\cdot), g_{1,m}(\cdot)$ and $g_{2,m}(\cdot)$ such that $p_m(h,e) \leq g_{0,m}(e)$, $|\frac{d}{dh}p_m(h,e)| \leq g_{1,m}(e)$, $|\frac{d^2}{dh^2}p_m(h,e)| \leq g_{2,m}(e)$ for all $(h,e)$. Furthermore, \linebreak $\sup_m \int |e|^r g_{i,m}(e) de < \infty $ for $i = 0,1,2$, and $0\leq r\leq 4$.
\end{enumerate}
\end{manualassumption}

\begin{manualassumption}{3.3} \label{interchange-assumption}
There exists a random variable $K(\xi(\theta))$ such that $\mathbb{E}_{\theta, \mathbb{P}_{\theta}}[K(\xi)] < \infty$, and the following holds in a probability 1 $(\mathbb{P}^N)$ subset of $\Theta$.
\begin{enumerate}
    \item $|h (x_2, \xi(\theta)) - h (x_1, \xi(\theta))| \leq K(\xi(\theta)) |x_2 - x_1|$ w.p.1 $(\mathbb{P}_{\theta})$ for all $x_1, x_2 \in \mathcal{X}$.
    \item The sample path gradient $d(x, \xi(\theta))$ exists w.p.1 $(\mathbb{P}_{\theta})$.
\end{enumerate}
\end{manualassumption}

\begin{manualassumption}{3.4} \label{hong-var-assumptions} \cite{Hong2009VaR}
    For any $x \in \mathcal{X}$, $H(x, \theta)$ has a continuous density $f(t; x)$ in a neighborhood of $t = v_{\alpha}(x)$, and $\partial_{x} F(t; x)$ exists and is continuous w.r.t. both $x$ and $t$ at $t = v_{\alpha}(x)$.
\end{manualassumption}

\begin{manualassumption}{3.5} \label{measure-assumption}
    Assume that there exists a family of measures $G_m(\cdot)$ and a number $\eta > 0$ such that for all $t \in \mathcal{B}_{\eta}(v_{\alpha}(x))$ and for all $\Delta y \subset (-\infty, \infty)$,
    \begin{equation*}
        |\nu(\Delta y, t) - \hat{\nu}^m (\Delta y, t)| \leq G_m(\Delta y) \text{ and } \int_{\mathbb{R}} |y| G_m(d y) \rightarrow 0 \text{ as } m \rightarrow \infty.
    \end{equation*}
\end{manualassumption}

\begin{manualassumption}{3.6} \label{supremum-assumption}
$\sup_{\theta} \mathbb{E}_{\mathbb{P}_{\theta}}[d(x; \xi)^2] < \infty$.
\end{manualassumption}

\begin{manualassumption}{3.10} \label{hong-assumption}
\cite{Hong2009CVaR}
\begin{enumerate}
    \item The VaR function $v_{\alpha}(x)$ is differentiable for any $x \in \mathcal{X}$.
    
    \item For any $x \in \mathcal{X}$, $P[H(x; \theta) = v_{\alpha}(x)] = 0$.
\end{enumerate}
\end{manualassumption}

\edit{
\section{Discussion on the assumptions}
In this paper, we present many technical assumptions that are needed for the results to hold. 
Here, we present a non-technical interpretation of each assumption to make it more intuitive to understand. 
Ignoring some pathological cases, we believe that most, if not all, assumptions are satisfied when $h(\cdot, \xi)$ is Lipschitz continuous for almost all $\xi$, the $\xi$ and $\theta$ are continuous random variables with light tails, and the variables $x, \xi, \theta$ are confined to compact (bounded) spaces.
\begin{itemize}
    \item Assumption 3.2: The random function $h(x, \xi)$ is smooth and has light tails. For any given $\theta$, the error function is a continuous random variable. 
    
    \item Assumption 3.3: The function $h(x, \xi)$ is Lipschitz continuous for a given $\xi$, and is differentiable almost everywhere.
    
    \item Assumption 3.4: In a neighborhood of its $\alpha$ quantile, the expected performance $H(x; \theta)$ is a continuous random variable (induced by $\theta$), and its distribution function is smooth in $x$.
    
    \item Assumption 3.5: This is a rather technical assumption that results from conditioning on measure zero events. It regulates the asymptotic behavior of the gradient observations, and requires the error of the gradient observations to converge to zero as the error of the function observation converges to zero.
    It is shown below that the assumption is satisfied for a general class of problems.
    
    \item Assumption 3.6: The gradient observations $d(x, \xi)$ have a finite second moment.
    
    \item Assumption 3.10: When we are interested in CVaR rather than VaR, this assumption weakens Assumption 3.4. It requires the VaR function to be differentiable, and the expected performance $H(x; \theta)$ to be a continuous random variable in a neighborhood of VaR.
    
\end{itemize}
}

\edit{
\subsection{A simple example}
The following example is modified from \cite{Hong2009VaR}. We proceed to show that the assumptions are satisfied, starting with the Assumption 3.5. 

\begin{example}
    Let $\theta_i \sim \mathcal{N}(0,1)$, and $h(x, \xi(\theta)) = x \theta_1 + \theta_2 + x \xi(\theta)$ where $\xi(\theta) \sim \mathcal{N}(0, \theta_1^2)$. Then, $H(x; \theta) = x \theta_1 + \theta_2$. 
    It follows that $D(x; \theta) = \theta_1$, $d(x, \xi(\theta)) = \theta_1 + \xi(\theta)$, $\hat{H}^m (x; \theta) = x \theta_1 + \theta_2 + x \frac{1}{m} \sum_{j=1}^m \xi_j(\theta)$ and $\hat{D}^m(x; \theta) = \theta_1 + \frac{1}{m} \sum_{j=1}^m \xi_j(\theta)$. 
\end{example}

To show that Assumption 3.5 holds, we have
\begin{align*}
    \nu((-\infty, y]; t) &= P(\theta_1 \leq y \mid x \theta_1 + \theta_2 = t) \\
    &= P \left( \theta_1 \leq y \mid \theta_1 = \frac{t - \theta_2}{x} \right) \\
    &= P(t - x y \leq \theta_2) \\
    & = 1 - \Phi(t - x y),
\end{align*}
where $\Phi(\cdot)$ is the CDF of $\mathcal{N}(0,1)$.
For the noisy counterpart, 
\begin{align*}
    \hat{\nu}^m((-\infty, y]; t) &= P \left(\theta_1 + \frac{1}{m} \sum_{j=1}^m \xi_j(\theta) \leq y \mid x \theta_1 + \theta_2 + x \frac{1}{m} \sum_{j=1}^m \xi_j(\theta) = t \right) \\
    &= P \left( \theta_1 + \frac{1}{m} \sum_{j=1}^m \xi_j(\theta) \leq y \mid \theta_1 + \frac{1}{m} \sum_{j=1}^m \xi_j(\theta) = \frac{t - \theta_2}{x} \right) \\
    &= P(t - x y \leq \theta_2) \\
    &= 1 - \Phi(t - x y).
\end{align*}
Thus, for this example, we have $\nu((-\infty, y]; t) = \hat{\nu}^m((-\infty]; t)$, and Assumption 3.5 is satisfied with $G_m(y) = 0$.
We now verify the remaining assumptions. 
\begin{itemize}
    \item Assumption 3.2: For the first part, $\tau_{\theta}^2 = x^2 \theta_1^2 + (x \theta_1 + \theta_2)^2 < \infty$, and $\tau^2 = 2 x^2 + 1 < \infty$. For the remaining parts, see the discussion in \cite{Gordy2010Nested} where it is implied that the assumption holds when the distribution of $h(x, \xi)$ is Gaussian.
    
    \item Assumption 3.3: The first part holds with $K(\xi(\theta)) = \theta_1 + \xi(\theta)$. For the second part, the derivative exists everywhere and is given by $d(x, \xi(\theta)) = \theta_1 + \xi(\theta)$.
    
    \item Assumption 3.4: $H(x; \theta)$ has a continuous density everywhere. Similarly $\partial_x F(t; x)$ exists and is continuous everywhere.
    
    \item Assumption 3.6: $\E_{\pth}[d(x, \xi)^2] = 2 \theta_1^2$. The supremum here is not finite, so the assumption is technically violated. However, this violation is purely technical and happens only due to unbounded domain of the normal distribution. We can get around this by truncating the domain $\Theta$ to some large interval, which is always done in practice due to limitations of  machine precision.
    
    \item Assumption 3.10: The VaR function is given by $v_{\alpha}(x) = z_{\alpha} \sqrt{x^2 + 1}$ where $z_{\alpha}$ is the $\alpha$ quantile of $\mathcal{N}(0, 1)$. It is differentiable for any $x \in \mathcal{X}$. The second part holds as $H(x; \theta)$ is a continuous random variable.
\end{itemize}
}

\edit{
\subsection{A general class of functions}
Here, we consider a class of functions of the form $h(x, \xi(\theta)) = H(x; \theta) + \xi(\theta)$ with $\xi(\theta)$ a mean zero finite variance random variable, and show that Assumption 3.5 holds for such problems. We will extend this to a more general class of functions below.
\begin{align*}
    \nu((-\infty, y]; t) &= P(D(x; \theta) \leq y \mid H(x; \theta) = t) \\
    \hat{\nu}^m((-\infty, y]; t) &= P(D(x; \theta) \leq y \mid H(x; \theta) + \Bar{\xi}(\theta) = t)
\end{align*}
If $H(x; \theta)$ is invertible with a continuous inverse function and $D(x; \theta)$ continuous in $\theta$, then the result follows as $\xi(\theta)$ converges to $0$, and the convergence rate is $\mathcal{O}(m^{-1/2})$.

For the case where the solution to $H(x; \theta) = t$ is not unique, we define $\vartheta(t) := \{\theta: H(x; \theta) = t \}$ as the restricted random variable that satisfies $H(x; \vartheta(t)) = t$. We require that $\vartheta(t)$ is sample path continuous, which follows if $H(x; \theta)$ is strictly monotone, i.e. has a non-zero gradient, in $\theta$ in a neighborhood of $t$. Then, 
\begin{align*}
    P(D(x; \theta) \leq y \mid H(x; \theta) + \Bar{\xi}(\theta) = t) &= \E_{\vartheta}[P(D(x; \vartheta(t - \Bar{\xi}(\vartheta))) \leq y)] \\
    & \xrightarrow{m \rightarrow \infty} \E_{\vartheta}[P(D(x; \vartheta(t)) \leq y)] \\
    &= P(D(x; \theta) \leq y \mid H(x; \theta) = t).
\end{align*}
Note that $\bar{\xi}(\theta) \xrightarrow{m \rightarrow \infty} 0$ uniformly in $\theta$ by Assumption 3.2. Thus, if $D(x; \theta)$ is Lipschitz continuous in $\theta$, the convergence is uniform in $y$, and the assumption follows.

\subsubsection{Extensions}
A simple extension is to functions of the form $h(x; \xi(\theta)) = H(x; \theta) + x \xi(\theta)$. The same analysis here holds, just with some extra terms which again disappear as $m \rightarrow \infty$. In this case, we have
\begin{align*}
    P(D(x; \theta) + \bar{\xi}(\theta) \leq y \mid H(x; \theta) + x \Bar{\xi}(\theta) = t) &= \E_{\vartheta}[P(D(x; \vartheta(t - x \Bar{\xi}(\vartheta))) + \bar{\xi}(\theta) \leq y)] \\
    & \xrightarrow{m \rightarrow \infty} \E_{\vartheta}[P(D(x; \vartheta(t)) \leq y)] \\
    &= P(D(x; \theta) \leq y \mid H(x; \theta) = t).
\end{align*}
This again follows by uniform convergence of $\bar{\xi}(\theta)$, and the Lipschitz continuity assumption on $D(x; \theta)$.

This line of argument extends to the functions of the form $h(x, \xi(\theta)) = H(x; \theta) + g(x) \xi(\theta)$ where $g(\cdot)$ is a differentiable function of $x$. Since the domain $\mathcal{X}$ is compact, the convergence is still uniform and the assumption holds. 

We note that the analysis presented here implicitly assumes that $\E_{\pth}[\xi] = 0$, which is not restrictive as we can always replace $\xi$ with $\xi - \E_{\pth}[\xi]$ to obtain a mean zero random variable. The limiting factor in this analysis is that we only considered error terms in which the degree of $\xi$ is one. When the error term includes higher order terms of $\xi$, the uniform convergence is no longer implied by Assumption 3.2, and this line of argument fails to hold. One can still check for uniform convergence, and use the same argument if it holds, or verify the assumption in some other way. We conclude this analysis by noting that when the error term is a first order function of $\xi$, the integral in Assumption 3.5 is $\mathcal{O}(m^{-1/2})$.
}

\edit{
\section{Details of numerical experiments}
In this section, we fill out the details of the numerical experiments that were left out from the paper due to space constraints.

\subsection{Details of the quadratic example}
In this subsection, we fill in the details of the quadratic example presented in Section 4.1 in the paper. We start by verifying the assumptions, continue with problem setup, and obtain an analytical solution to the BRO optimization problem.

This example is modified from the simple example by \cite{Hong2009VaR} presented above. Assuming a posterior distribution of the form $\theta_i \sim \mathcal{N}(\mu_i, \sigma_i^2)$, the assumptions can be verified in a similar way. Here, we briefly highlight the differences.
\begin{itemize}
    \item Assumption 3.2: $\tau_{\theta}^2 = x^2 \theta_1^2 + (x \theta_1 + x^2 \theta_2)^2 < \infty$ and $\tau^2 = 2 x^2 (\sigma_1^2 + \mu_1^2) + 2 x^3 \mu_1 \mu_2 + x^4 (\sigma_2^2 + \mu_2^2) < \infty$.
    
    \item Assumption 3.3: We have $K(\xi(\theta)) = \sup_{x \in \mathcal{X}} \theta_1 + 2x \theta_2 + \xi(\theta)$ which is integrable since $\mathcal{X}$ is compact. The derivative is given above and exists everywhere.
    
    \item Assumption 3.4: Same as before.
    
    \item Assumption 3.5: The assumption is satisfied with $G_m(\cdot) = 0$ as shown below.
    \begin{align*}
        \nu((-\infty, y]; t) &= P(D(x; \theta) \leq y \mid H(x; \theta) = t) \\
        &= P(\theta_1 + 2x \theta_2 \leq y \mid x \theta_1 + x^2 \theta_2 = t) \\
        &= P \left( \theta_1 + 2x \theta_2 \leq y \mid \theta_1 = \frac{t - x^2 \theta_2}{x} \right) \\
        &= P(t + x^2 \theta_2 \leq xy) \\
        &= P \left(\theta_2 \leq \frac{xy - t}{x^2} \right)
    \end{align*}
    \begin{align*}
        \hat{\nu}^m((-\infty, y]; t) &= P(\hat{D}^m(x; \theta) \leq y \mid \hat{H}^m(x; \theta) = t) \\
        &= P(\theta_1 + 2x \theta_2 + \bar{\xi}(\theta) \leq y \mid x \theta_1 + x^2 \theta_2 + x \bar{\xi}(\theta) = t ) \\
        &= P \left( \theta_1 + 2x \theta_2 + \bar{\xi}(\theta) \leq y \mid \theta_1 + \bar{\xi}(\theta) = \frac{t - x^2 \theta_2}{x} \right) \\
        &= P\left( \theta_2 \leq \frac{xy - t}{x^2} \right)
    \end{align*}
    Thus, $|\nu((-\infty, y]; t) - \hat{\nu}^m((-\infty, y]; t)| = 0$, and $G_m(\cdot) = 0$ satisfies the assumption.
    
    \item Assumption 3.6: We again run into same technical violation, which is resolved by bounding the domain $\Theta$.
    
    \item Assumption 3.10: The expression for the VaR is given below and is continuously differentiable almost everywhere. The second part holds as $H(x; \theta)$ is a continuous random variable.
\end{itemize}

With the given posterior of the form $\theta_i \sim \mathcal{N}(\mu_i, \sigma_i^2)$, $H(x; \theta)$ is a normal random variable, and the VaR and CVaR objectives are given by 
\begin{equation*}
    v_{\alpha}(x) = \mu_H + z_{\alpha} \sigma_H \text{ and } c_{\alpha}(x) = \mu_H + \frac{\phi(z_{\alpha})}{1-\alpha} \sigma_H,
\end{equation*}
where $z_{\alpha}$ and $\phi(\cdot)$ are the $\alpha$ quantile and the PDF of the standard normal distribution respectively, and $\mu_H, \sigma_H$ are the mean and standard deviation of $H(x; \theta)$ given by
\begin{equation*}
    \mu_H = x \mu_1 + x^2 \mu_2 \text{ and } \sigma_H = \sqrt{x^2 \sigma_1^2 x^4 \sigma_2^2 }.
\end{equation*}
The gradients of VaR and CVaR can be computed in a similar manner, and are given by $v'_{\alpha}(x) = \mu'_H + z_{\alpha} \sigma'_H$ and $c'_{\alpha}(x) = \mu'_H + \frac{\phi(z_{\alpha})}{1-\alpha} \sigma'_H$, where $\mu'_H$ and $\sigma'_H$ denote the gradients of $\mu_H$ and $\sigma_H$ respectively.

For a given $\theta = (\theta_1, \theta_2)$ with $\theta_2 > 0$, the minimizer is given by $x^*(\theta) = -\frac{\theta_1}{2 \theta_2}$. Once $(\mu_1, \mu_2, \sigma_1, \sigma_2)$ are specified, the minimizers of BRO objectives can also be computed analytically. 

We consider the case where $\theta_1, \theta_2$ are estimated from samples $\{\xi_i^j\}_{j = 1, \ldots, N} \stackrel{iid}{\sim} \mathcal{N}(\theta_i, \tilde{\sigma}_i^2)$ where $\tilde{\sigma}_i$ is the known variance. Then, using a degenerate normal prior on $\theta$, the posterior distribution is given by $\theta_i \sim \mathcal{N}(\bar{\xi_i}, \frac{\tilde{\sigma}_i^2}{N})$. In this numerical example, we suppose that the posterior is given by $\theta_1 \sim \mathcal{N}(-15, 16)$ and $\theta_2 \sim \mathcal{N}(10, 4)$, and consider the BRO objective with risk measure CVaR at risk level $\alpha = 0.75$. With the given parameters, the optimal solution is found at $x^* = 0.474775$ with the corresponding BRO CVaR objective value of $-2.38647$.

To keep things simple, we use a fixed budget sequence of $n_t = n = 100$ and $m_t = m = n / 5$. The gradient-based algorithms use the estimators developed in this paper, and the gradient-free alternatives use the nested estimators of \cite{Zhu2017input-uncertainty} to estimate the objective value. The benchmark algorithms we consider are originally developed for deterministic optimization. Thus, in addition to stochastic evaluations of the objective, we also consider the Sample Average Approximation (SAA, \citealp{Kim2015SAA}) which converts the stochastic optimization problem into an approximate deterministic optimization problem. This is done by fixing a random draw of $\theta_i$ and $\xi_j(\theta_i)$ before the optimization starts, and using this fixed set of samples to calculate the value of the estimators. For comparison, the stochastic evaluations draw a new set of random variables $\theta_i$ and $\xi_j(\theta_i)$ for each evaluation of the estimators. A different set of $\theta_i$ and $\xi_j(\theta_i)$ is used for each replication of SAA.

\subsection{Convergence analysis of MCMC output}
It is known that MCMC methods converge to a steady state distribution, however, detecting when this convergence occurs, without access to the distribution, is an open question. In this paper, we use an improvised Wasserstein distance analysis to perform an empirical convergence analysis. 

The idea behind this analysis is as follows. We treat subsets of the MCMC chain as empirical distributions drawn from the chain's distribution at a given time. The subsets should be sufficiently large so that the correlation between the samples can be ignored. If the chain has not converged to the steady state distribution, the underlying distribution is actively changing and this should show up as a distance between the two empirical distributions. As the convergence occurs, we can expect the distance between subsequent empirical distributions to get smaller, and stabilize as we converge to the stationary distribution. Note that, even when the empirical distributions are drawn from the stationary distribution, the distance between the empirical distributions will be non-zero as these are essentially two random draws from the same underlying distribution. 

However, once we have reached the stationary distribution, the distance between any two arbitrary subsets, not necessarily subsequent, should be roughly the same, apart from the random noise. We observe this behavior in the empirical posterior distributions used in the experiments. The Wasserstein distance between subsets of size $10^5$ are observed to be about $10^{-3}$ ($\pm$ noise) regardless of the ordering of the subsets. Thus, we conclude that the underlying distribution is not actively changing and the empirical distribution can be treated as coming from the steady state, i.e. the true posterior, distribution.
}

\section{Results and proofs}
Since the results are given for a fixed $x$, we drop the dependence on $x$ to simplify notation in the proofs. Moreover, $\hat{H}^m_i$ and $H_i$ are commonly used in place of $\hat{H}^m(x; \theta_i)$ and $H(x; \theta_i)$ respectively. We do the same with $\hat{D}^m_i$ and $D_i$ as well.


\begin{manualproposition}{3.7} \label{var-unbiased-proposition}
    Suppose that Assumptions \ref{zhu-assumption}, \ref{interchange-assumption}, \ref{hong-var-assumptions}, \ref{measure-assumption} and \ref{supremum-assumption} hold. Then $\mathbb{E}_{\theta, \mathbb{P}_{\theta}}[\varphi^{n,m}_{\alpha}(x)] - v'_{\alpha}(x) \rightarrow 0$ as $n,m \rightarrow \infty$.
    
    Moreover, if in addition the integral in Assumption \ref{measure-assumption} is $\mathcal{O}(m^{-1/2})$, $g(t; x)$ is differentiable w.r.t. $t$ at $t = v_{\alpha}(x)$, and the budget sequence is such that $n = \Theta(m)$, then the bias is $\mathbb{E}_{\theta, \mathbb{P}_{\theta}}[\varphi^{n,m}_{\alpha}(x)] - v'_{\alpha}(x) = \mathcal{O}(n^{-1/2})$.
\end{manualproposition}

\begin{proof}
We will show that $\lim_{n, m \rightarrow \infty} \mathbb{E}_{\theta, \mathbb{P}_{\theta}}[\varphi^{n,m}_{\alpha}] = \mathbb{E}_{\theta}[\partial_x H(x; \theta) \mid H(x; \theta) = v_{\alpha}]$. Let $F_v$ denote the CDF of $\hat{v}^{n,m}_{\alpha}$, and let $\hat{g}^m(t; x) = \mathbb{E}_{\theta, \mathbb{P}_{\theta}}[\partial_x \hat{H}^m (x; \theta) \mid \hat{H}^m (x; \theta) = t]$. Using an argument from the proof of Theorem 3 in \cite{Hong2009VaR},
\begin{align*}
    \mathbb{E}_{\theta, \mathbb{P}_{\theta}}[\varphi^{n,m}_{\alpha}] &= \mathbb{E}_{\theta, \mathbb{P}_{\theta}}[\hat{D}^m (\hat{\theta}^m_{(\lceil \alpha n \rceil)})]\\
    &=\int \mathbb{E}_{\theta, \mathbb{P}_{\theta}}[\hat{D}^m (\hat{\theta}^m_{(\lceil \alpha n \rceil)}) \mid \hat{H}^m (\hat{\theta}^m_{(\lceil \alpha n \rceil)}) = t] dF_v(t) \tag{$\dagger$}\\
    &=\int \mathbb{E}_{\theta, \mathbb{P}_{\theta}}[\hat{D}^m (\theta) \mid \hat{H}^m(\theta) = t] dF_v(t) \tag{$\dagger\dagger$} \\
    &= \int \hat{g}^m (t;x) dF_v(t) \\
    &= \mathbb{E}_{\theta, \mathbb{P}_{\theta}}[\hat{g}^m (\hat{v}^{n,m}_{\alpha}; x)],
\end{align*}
where $(\dagger\dagger)$ follows from the fact that $\theta_1, \theta_2, \ldots, \theta_n$ are i.i.d., and thus $\{\hat{H}^m(\theta_i)\}_{i=1}^n$ are i.i.d.; conditioning on $\{\hat{H}^m(\hat{\theta}^m_{(\lceil \alpha n \rceil)}) = t\}$ is the same as conditioning on 
\begin{enumerate}
    \item there exists some $i^*$ such that $\hat{H}^m(\theta_{i^*}) = t$;
    
    \item for $i \neq i^*$, there are $\alpha n -1$ values of $\hat{H}^m (\theta_i)$ which are $< t$, and the rest are $>t$.
\end{enumerate}
Moreover, due to independence and symmetry, $(\dagger)$ is the same as (without loss of generality) $\mathbb{E}_{\theta, \mathbb{P}_{\theta}}[ \hat{D}^m(\theta_1) \mid \hat{H}^m(\theta_1) = t]$ where $\theta_1 \sim \mathbb{P}^N$ a random variable and not a fixed realization. We refer the reader to \cite{Hong2009VaR} for more details.

Next, we show that $\hat{g}^m(\cdot; x) \rightarrow g(\cdot; x)$ uniformly on $\mathcal{B}_\eta(v_{\alpha})$ as $m \rightarrow \infty$. Note that we can write,
\begin{equation*}
    g(t;x) = \mathbb{E}_{\theta}[\partial_x H(x; \theta) \mid H(x; \theta) = t] = \int_{\mathbb{R}} y \nu(d y, t);
\end{equation*}
and similarly,
\begin{equation*}
    \hat{g}^m(t; x) = \int_{\mathbb{R}} y \hat{\nu}^m (d y, t).
\end{equation*}
Then, for all $t \in \mathcal{B}_{\eta}(v_{\alpha})$,
\begin{align*}
    |\hat{g}^m (t;x) - g(t;x)| &= \left| \int_{\mathbb{R}} y \left( \hat{\nu}^m(d y,t) - \nu
    (d y,t) \right) \right| \\
    &\leq \int_{\mathbb{R}} |y| \left| \hat{\nu}^m(d y,t) - \nu
    (d y,t) \right| \\
    &\leq \int_\mathbb{R} |y| G_m(d y) \rightarrow 0 \quad \text{ as } m \rightarrow \infty,
\end{align*}
where the first inequality follows from the triangle inequality and the definition of Lebesgue integral, and the last inequality follows from Assumption \ref{measure-assumption}. Since $G_m(\cdot)$ does not depend on $t$, the convergence is uniform in $t$. We now claim that $\hat{g}^m (\hat{v}^{n,m}_{\alpha}; x) \rightarrow g(v_{\alpha}; x)$ a.s. as $n,m \rightarrow \infty$. To see this, take a sample path on which $\hat{v}^{n,m}_{\alpha} \rightarrow v_{\alpha}$ as $n, m \rightarrow \infty$, which occurs a.s. due to Assumption \ref{zhu-assumption}. Therefore, it suffices to show convergence on this sample path. Take any $\epsilon>0$, and notice the following.
\begin{itemize}
    \item Since $g(\cdot; x)$ is continuous by Assumption \ref{hong-var-assumptions}, there exists $\delta_1 > 0$  such that \begin{equation*}
        |g(t; x) - g(v_{\alpha}; x)| < \epsilon /2, \quad \forall t \in \mathcal{B}_{\delta_1}(v_{\alpha}).
    \end{equation*}
    
    \item Take $\delta_2 := \min \{ \delta, \eta \}$. Then, there exists $N_1, M_1 \in \mathbb{Z}^+$ such that 
    \begin{equation*}
        |\hat{v}^{n,m}_{\alpha} - v_{\alpha}| < \delta_2, \quad \forall n \geq N_1, m \geq M_1.
    \end{equation*}
    
    \item Furthermore, due to uniform convergence, there exists $M_2 \in \mathbb{Z}^+$ such that 
    $\forall m \geq M_2$, 
    \begin{equation*}
        |\hat{g}^m(t; x) - g(t;x)| < \epsilon / 2, \quad \forall t \in \mathcal{B}_{\delta_2}(v_{\alpha}).
    \end{equation*}
\end{itemize}
Combining the above, we get
\begin{align*}
    &|\hat{g}^m (\hat{v}^{n,m}_{\alpha}; x) - g(v_{\alpha};x)| \\
    \leq& |\hat{g}^m (\hat{v}^{n,m}_{\alpha}; x) - g(\hat{v}^{n,m}_{\alpha}; x)| + |g(\hat{v}^{n,m}_{\alpha}; x) - g(v_{\alpha}; x)| \\
    \leq& \frac{\epsilon}{2} + \frac{\epsilon}{2} = \epsilon,
\end{align*}
for all $n \geq N_1, m \geq \max \{ M_1, M_2 \}$. Thus, 
\begin{equation} \label{eq:as-conv}
    \hat{g}^m(\hat{v}^{n,m}_{\alpha}; x) \rightarrow g(v_{\alpha}; x) \quad \text{ a.s. as }  n,m \rightarrow \infty. 
\end{equation}
Finally, since
\begin{align*}
    \mathbb{E}_{\theta, \mathbb{P}_{\theta}}[\hat{g}^m (\hat{v}^{n,m}_{\alpha}; x)^2] &= \int \mathbb{E}_{\theta, \mathbb{P}_{\theta}}[\partial_x \hat{H}^m(x; \theta) \mid \hat{H}^m(x; \theta) = t] dF_v(t) \\
    &\leq \int \mathbb{E}_{\theta, \mathbb{P}_{\theta}} \left\{[\partial_x \hat{H}^m (x; \theta)]^2 \mid \hat{H}^m (x; \theta) = t \right\} dF_v(t) \tag{Jensen's inequality} \\
    &= \mathbb{E}_{\theta, \mathbb{P}_{\theta}}\left\{ [\partial_x \hat{H}^m(x; \hat{\theta}^m_{(\lceil \alpha n \rceil)})]^2 \right\} \tag{Similar to $(\dagger\dagger)$} \\
    &\leq \sup_{\theta} \mathbb{E}_{\mathbb{P}_{\theta}}[d(x; \xi)^2] < \infty, \tag{By Assumption \ref{supremum-assumption}}
\end{align*}
we get that $\hat{g}^m(\hat{v}^{n,m}_{\alpha}; x)$ is uniformly integrable. This together with (\ref{eq:as-conv}) yields that
\begin{equation*}
\lim_{n, m \rightarrow \infty} \mathbb{E}_{\theta, \mathbb{P}_{\theta}}[\varphi^{n,m}_{\alpha}]  = g(v_{\alpha};x) = v'_{\alpha}.
\end{equation*}
This completes the first part of the proof. For the second part of the proposition, given that the integral is $\mathcal{O}(m^{-1/2})$, we have that $\hat{g}^m(t; x) - g(t; x) = \mathcal{O}(m^{-1/2})$ uniformly for all $t \in \mathcal{B}_{\eta}(v_{\alpha})$. Then, for $n, m$ large enough (so that $\hat{v}^{n,m}_{\alpha} \in \mathcal{B}_{\eta}(v_{\alpha})$),
\begin{align*}
    |\hat{g}^m(\hat{v}^{n,m}_{\alpha}; x) - g(v_{\alpha}; x)| &\leq |g(\hat{v}^{n,m}_{\alpha}; x) - g(v_{\alpha}; x)| + \mathcal{O}(m^{-1/2}).
\end{align*}
For any fixed $\hat{v}^{n,m}_{\alpha}$, using Taylor's theorem,
\begin{align*}
    |g(\hat{v}^{n,m}_{\alpha}; x) - g(v_{\alpha}; x)| = |\partial_t g(v_{\alpha}; x) (\hat{v}^{n,m}_{\alpha} - v_{\alpha}) + o(\hat{v}^{n,m}_{\alpha} - v_{\alpha})|
\end{align*}
Putting it all together, we have
\begin{align*}
    \mathbb{E}_{\theta, \mathbb{P}_{\theta}}[\varphi^{n,m}_{\alpha}] - v'_{\alpha} &= \mathbb{E}_{\theta, \mathbb{P}_{\theta}}[\hat{g}^m(\hat{v}^{n,m}_{\alpha}; x)] - g(v_{\alpha}; x) \\
    &\leq \mathbb{E}_{\theta, \mathbb{P}_{\theta}}[|g(\hat{v}^{n,m}_{\alpha}; x) - g(v_{\alpha}; x)| ] + \mathcal{O}(m^{-1/2}) \\
    &= \mathbb{E}_{\hat{v}^{n,m}_{\alpha}} [ | \partial_t g(v_{\alpha}; x) (\hat{v}^{n,m}_{\alpha} - v_{\alpha}) + o(\hat{v}^{n,m}_{\alpha} - v_{\alpha}) | ] + \mathcal{O}(m^{-1/2}).
\end{align*}
Theorem 3.6 of \cite{Zhu2017input-uncertainty} shows that under Assumption \ref{zhu-assumption}, $n = o(m^2)$ is a sufficient and necessary condition for 
\begin{equation*}
    \sqrt{n} (\hat{v}^{n,m}_{\alpha} - v_{\alpha}) \Rightarrow \mathcal{N}(0, \sigma_v)
\end{equation*}
where $\sigma_v = \frac{\alpha (1-\alpha)}{f^{v_{\alpha}}}$. Therefore, the term inside the absolute value (when scaled by $\sqrt{n}$) converges to a mean zero normal random variable, and the expectation is $\mathcal{O}(n^{-1/2})$. Putting it together with $n = \Theta(m)$, we get that the bias is $\mathcal{O}(n^{-1/2})$.
\end{proof}



\begin{manualtheorem}{3.8} \label{var-consistency-theorem}
    Suppose that Assumptions \ref{zhu-assumption}, \ref{interchange-assumption}, \ref{hong-var-assumptions}, \ref{measure-assumption} and \ref{supremum-assumption} hold, then 
    \begin{equation*}
        \bar{\varphi}^{n,m,k}_{\alpha}(x) \xrightarrow{P} v'_{\alpha}(x) \text{ as } n, m, k \rightarrow \infty,
    \end{equation*}
    where $\xrightarrow{P}$ denotes convergence in probability.
\end{manualtheorem}

\begin{proof}
By Proposition \ref{var-unbiased-proposition}, we have $\mathbb{E}_{\theta, \mathbb{P}_{\theta}}[\varphi^{n,m}_{\alpha}] \rightarrow v'_{\alpha}$ as $n, m \rightarrow \infty$. For any $\epsilon > 0$, by Chebyshev's inequality,
\begin{align*}
    P \left( \left| \bar{\varphi}^{n,m,k}_{\alpha} - \mathbb{E}_{\theta, \mathbb{P}_{\theta}} [\varphi^{n,m}_{\alpha}] \right| > \epsilon \right) &\leq \frac{Var(\varphi^{n,m}_{\alpha})}{k \epsilon^2} 
    \leq \frac{\mathbb{E}_{\theta, \mathbb{P}_{\theta}}[(\varphi^{n,m}_{\alpha})^2]}{k \epsilon^2} 
    \leq \frac{\sup_{\theta} \mathbb{E}_{\mathbb{P}_{\theta}}[d(\xi)^2]}{k \epsilon^2}.
\end{align*}
Therefore, $\bar{\varphi}^{n,m,k}_{\alpha} \xrightarrow{P} \mathbb{E}[\varphi^{n,m}_{\alpha}]$ as $k \rightarrow \infty$ uniformly for each $n, m$. Combining with the result of Proposition \ref{var-unbiased-proposition}, $\bar{\varphi}^{n,m,k}_{\alpha} \xrightarrow{P} v'_{\alpha}$ as $n, m, k \rightarrow \infty$.
\end{proof}


\begin{manualproposition}{3.11} \label{indicator-proposition}
Suppose Assumption \ref{zhu-assumption} holds and $P(H(x; \theta) = v_{\alpha}(x)) = 0$. Then
$$ \frac{1}{n} \sum_{i=1}^n | \mathds{1}_{(\hat{H}^m (x; \theta_i) \geq \hat{v}^{n,m}_{\alpha}(x))} - \mathds{1}_{(H(x; \theta_i) \geq v_{\alpha}(x))} | \rightarrow 0 \ w.p.1 \text{ as } n,m \rightarrow \infty .$$
\end{manualproposition}

\begin{proof}

\begin{align*}
\frac{1}{n} \sum_{i=1}^n | \mathds{1}_{(\hat{H}^m_i \geq \hat{v}^{n,m}_{\alpha})} - \mathds{1}_{(H_i \geq v_{\alpha})} |
&= \frac{1}{n} \sum_{i=1}^n | \mathds{1}_{(\hat{H}^m_i \geq \hat{v}^{n,m}_{\alpha})} - \mathds{1}_{(H_i \geq \hat{v}^{n,m}_{\alpha})} + \mathds{1}_{(H_i \geq \hat{v}^{n,m}_{\alpha})} - \mathds{1}_{(H_i \geq v_{\alpha})} | \\
&\leq \underbrace{\frac{1}{n} \sum_{i=1}^n | \mathds{1}_{(\hat{H}^m_i \geq \hat{v}^{n,m}_{\alpha})} - \mathds{1}_{(H_i \geq \hat{v}^{n,m}_{\alpha})} |}_{(P1)} + \underbrace{ \frac{1}{n} \sum_{i=1}^n | \mathds{1}_{(H_i \geq \hat{v}^{n,m}_{\alpha})} - \mathds{1}_{(H_i \geq v_{\alpha})} | }_{(P2)} 
\end{align*}
We will now proceed to bound each term individually and show that the sum is bounded above by zero in the limit.
\begin{align*}
(P2) &= ( \mathds{1}_{( \hat{v}^{n,m}_{\alpha} \leq v_{\alpha} )} - \mathds{1}_{( \hat{v}^{n,m}_{\alpha} > v_{\alpha} )} ) \frac{1}{n} \sum_{i=1}^n ( \mathds{1}_{(H_i \geq \hat{v}^{n,m}_{\alpha})} - \mathds{1}_{(H_i \geq v_{\alpha})} ) \\
&= \mathds{1}_{( \hat{v}^{n,m}_{\alpha} \leq v_{\alpha} )} \underbrace{ \frac{1}{n} \sum_{i=1}^n ( \mathds{1}_{(H_i \geq \hat{v}^{n,m}_{\alpha})} - \mathds{1}_{(H_i \geq v_{\alpha})} ) }_{(P3)} + \mathds{1}_{( \hat{v}^{n,m}_{\alpha} > v_{\alpha} )} \underbrace{ \frac{1}{n} \sum_{i=1}^n ( \mathds{1}_{(H_i \geq v_{\alpha})} - \mathds{1}_{(H_i \geq \hat{v}^{n,m}_{\alpha})} ) }_{(P4)}
\end{align*}
Note that both $(P3)$ and $(P4)$ are non-negative as long as the accompanying indicator is one.
\begin{align*}
(P3) &= \mathds{1}_{( \hat{v}^{n,m}_{\alpha} < v_{\alpha} - \epsilon )} (P3) + \mathds{1}_{( \hat{v}^{n,m}_{\alpha} \geq v_{\alpha} - \epsilon )} (P3) \\
&\leq \mathds{1}_{( \hat{v}^{n,m}_{\alpha} < v_{\alpha} - \epsilon )} + \mathds{1}_{( \hat{v}^{n,m}_{\alpha} \geq v_{\alpha} - \epsilon )} \frac{1}{n} \sum_{i=1}^n ( \mathds{1}_{( H_i \geq v_{\alpha} - \epsilon )} - \mathds{1}_{( H_i \geq v_{\alpha})} ) \\
&\leq \mathds{1}_{( \hat{v}^{n,m}_{\alpha} < v_{\alpha} - \epsilon )} + \frac{1}{n} \sum_{i=1}^n  \mathds{1}_{( v_{\alpha} > H_i \geq v_{\alpha} - \epsilon )}
\end{align*}
Similarly,
\begin{align*}
(P4) &= \mathds{1}_{( \hat{v}^{n,m}_{\alpha} \leq v_{\alpha} + \epsilon )} (P4) + \mathds{1}_{( \hat{v}^{n,m}_{\alpha} > v_{\alpha} + \epsilon )} (P4) \\
&\leq \mathds{1}_{( \hat{v}^{n,m}_{\alpha} > v_{\alpha} + \epsilon )} + \mathds{1}_{( \hat{v}^{n,m}_{\alpha} \leq v_{\alpha} + \epsilon )} \frac{1}{n} \sum_{i=1}^n ( \mathds{1}_{( H_i \geq v_{\alpha} )} - \mathds{1}_{( H_i \geq v_{\alpha} + \epsilon)} ) \\
&\leq \mathds{1}_{( \hat{v}^{n,m}_{\alpha} > v_{\alpha} + \epsilon )} + \frac{1}{n} \sum_{i=1}^n \mathds{1}_{( v_{\alpha} + \epsilon > H_i \geq v_{\alpha} )} 
\end{align*}
Putting them together, we get
\begin{align*}
(P2) &\leq  \mathds{1}_{( \hat{v}^{n,m}_{\alpha} < v_{\alpha} - \epsilon )} + \mathds{1}_{( \hat{v}^{n,m}_{\alpha} > v_{\alpha} + \epsilon )} + \sum_{i=1}^n ( \mathds{1}_{( v_{\alpha} > H_i \geq v_{\alpha} - \epsilon )} + \mathds{1}_{( v_{\alpha} + \epsilon > H_i \geq v_{\alpha} )} ) \\
&= (1 - \mathds{1}_{( v_{\alpha} + \epsilon \geq \hat{v}^{n,m}_{\alpha} \geq v_{\alpha} - \epsilon )} ) + \frac{1}{n} \sum_{i=1}^n \mathds{1}_{( v_{\alpha} + \epsilon > H_i \geq v_{\alpha} - \epsilon )}
\end{align*}
We can now look at the other term.
\begin{align*}
(P1) &= \frac{1}{n} \sum_{i=1}^n | \mathds{1}_{( \hat{H}^m_i \geq \hat{v}^{n,m}_{\alpha} )} - \mathds{1}_{( H_i \geq \hat{v}^{n,m}_{\alpha} )} | \\
&= \frac{1}{n} \sum_{i=1}^n ( \mathds{1}_{( \hat{H}^m_i \geq \hat{v}^{n,m}_{\alpha} > H_i )} + \mathds{1}_{( H_i \geq \hat{v}^{n,m}_{\alpha} > \hat{H}^m_i )} ) \\
&= \frac{1}{n} \sum_{i=1}^n \left[ \mathds{1}_{( \hat{H}^m_i \leq H_i + \epsilon )} \mathds{1}_{( \hat{H}^m_i \geq \hat{v}^{n,m}_{\alpha} > H_i )} + \mathds{1}_{( \hat{H}^m_i > H_i + \epsilon )} \mathds{1}_{( \hat{H}^m_i \geq \hat{v}^{n,m}_{\alpha} > H_i )} \right. \\
&\qquad \qquad \left. + \mathds{1}_{( \hat{H}^m_i \geq H_i - \epsilon )} \mathds{1}_{( H_i \geq \hat{v}^{n,m}_{\alpha} > \hat{H}^m_i )} + \mathds{1}_{( \hat{H}^m_i < H_i - \epsilon )} \mathds{1}_{( H_i \geq \hat{v}^{n,m}_{\alpha} > \hat{H}^m_i )} \right] \\
&\leq \underbrace{ \frac{1}{n} \sum_{i=1}^n \left[ \mathds{1}_{( \hat{H}^m_i \leq H_i + \epsilon )} \mathds{1}_{( H_i + \epsilon \geq \hat{v}^{n,m}_{\alpha} > H_i )} + \mathds{1}_{( \hat{H}^m_i \geq H_i - \epsilon )} \mathds{1}_{( H_i \geq \hat{v}^{n,m}_{\alpha} > H_i - \epsilon )} \right] }_{(P5)} \\
&\qquad \qquad + \frac{1}{n} \sum_{i=1}^n ( \mathds{1}_{( \hat{H}^m_i > H_i + \epsilon )} + \mathds{1}_{( \hat{H}^m_i < H_i - \epsilon )} )
\end{align*}

\begin{align*}
(P5) &\leq \frac{1}{n} \sum_{i=1}^n ( \mathds{1}_{( H_i + \epsilon \geq \hat{v}^{n,m}_{\alpha} > H_i )} + \mathds{1}_{( H_i \geq \hat{v}^{n,m}_{\alpha} > H_i - \epsilon )} ) \\
&= \frac{1}{n} \sum_{i=1}^n \mathds{1}_{( H_i + \epsilon \geq \hat{v}^{n,m}_{\alpha} > H_i - \epsilon )} = \frac{1}{n} \sum_{i=1}^n \mathds{1}_{( H_i + \epsilon \geq \hat{v}^{n,m}_{\alpha} )} \mathds{1}_{(  \hat{v}^{n,m}_{\alpha} > H_i - \epsilon )} \\
&= \frac{1}{n} \sum_{i=1}^n \left[ \mathds{1}_{( v_{\alpha} + \epsilon \geq \hat{v}^{n,m}_{\alpha} \geq v_{\alpha} - \epsilon )} \mathds{1}_{( H_i + \epsilon \geq \hat{v}^{n,m}_{\alpha} )} \mathds{1}_{(  \hat{v}^{n,m}_{\alpha} > H_i - \epsilon )} \right. \\
&\qquad \qquad \left. + (1 - \mathds{1}_{( v_{\alpha} + \epsilon \geq \hat{v}^{n,m}_{\alpha} \geq v_{\alpha} - \epsilon )}) \mathds{1}_{( H_i + \epsilon \geq \hat{v}^{n,m}_{\alpha} )} \mathds{1}_{(  \hat{v}^{n,m}_{\alpha} > H_i - \epsilon )} \right] \\
&\leq \frac{1}{n} \sum_{i=1}^n  \mathds{1}_{( v_{\alpha} + \epsilon \geq \hat{v}^{n,m}_{\alpha} \geq v_{\alpha} - \epsilon )} \mathds{1}_{( H_i + \epsilon \geq \hat{v}^{n,m}_{\alpha} )} \mathds{1}_{(  \hat{v}^{n,m}_{\alpha} > H_i - \epsilon )} + (1 - \mathds{1}_{( v_{\alpha} + \epsilon \geq \hat{v}^{n,m}_{\alpha} \geq v_{\alpha} - \epsilon )})   \\
&\leq \frac{1}{n} \sum_{i=1}^n  \mathds{1}_{( H_i + 2 \epsilon \geq v_{\alpha} > H_i - 2 \epsilon )}  + (1 - \mathds{1}_{( v_{\alpha} + \epsilon \geq \hat{v}^{n,m}_{\alpha} \geq v_{\alpha} - \epsilon )}) 
\end{align*}
Putting them together,
\begin{align*}
(P1) &\leq (1 - \mathds{1}_{( v_{\alpha} + \epsilon \geq \hat{v}^{n,m}_{\alpha} \geq v_{\alpha} - \epsilon )}) + \frac{1}{n} \sum_{i=1}^n  \mathds{1}_{( H_i + 2 \epsilon \geq v_{\alpha} > H_i - 2 \epsilon )} + \frac{1}{n} \sum_{i=1}^n ( \mathds{1}_{( \hat{H}^m_i > H_i + \epsilon )} + \mathds{1}_{( \hat{H}^m_i < H_i - \epsilon )}) \\
&= (1 - \mathds{1}_{( v_{\alpha} + \epsilon \geq \hat{v}^{n,m}_{\alpha} \geq v_{\alpha} - \epsilon )}) + \frac{1}{n} \sum_{i=1}^n  \mathds{1}_{( H_i + 2 \epsilon \geq v_{\alpha} > H_i - 2 \epsilon )} + \frac{1}{n} \sum_{i=1}^n ( 1 - \mathds{1}_{( H_i + \epsilon \geq \hat{H}^m_i \geq H_i - \epsilon )} )
\end{align*}
Thus, we have 
\begin{align*}
\frac{1}{n} \sum_{i=1}^n | \mathds{1}_{(\hat{H}^m_i \geq \hat{v}^{n,m}_{\alpha})} - \mathds{1}_{(H_i \geq v_{\alpha})} | &\leq 2 (1 - \mathds{1}_{( v_{\alpha} + \epsilon \geq \hat{v}^{n,m}_{\alpha} \geq v_{\alpha} - \epsilon )}) + \frac{1}{n} \sum_{i=1}^n  \mathds{1}_{( H_i + 2 \epsilon \geq v_{\alpha} > H_i - 2 \epsilon )} \\
&\qquad \qquad + \frac{1}{n} \sum_{i=1}^n ( 1 - \mathds{1}_{( H_i + \epsilon \geq \hat{H}^m_i \geq H_i - \epsilon )} ) + \frac{1}{n} \sum_{i=1}^n \mathds{1}_{( v_{\alpha} + \epsilon > H_i \geq v_{\alpha} - \epsilon )}
\end{align*}
By strong consistency of $\hat{v}^{n,m}_{\alpha}$, the first term goes to $0$ w.p.1. The second term is equivalent to $P(H + 2\epsilon \geq v_{\alpha} > H - 2 \epsilon)$ (w.p.1) in the limit which goes to $0$ as $\epsilon \rightarrow 0$ since $P(H = v_{\alpha}) = 0$. The last term goes to $0$ w.p.1 by a similar argument.
Let's focus on the third term. We can rewrite it as
\begin{align*}
    (*) = \frac{1}{n} \sum_{i=1}^n ( 1 - \mathds{1}_{( H_i + \epsilon \geq \hat{H}^m_i \geq H_i - \epsilon )} ) &= \frac{1}{n} \sum_{i=1}^n \mathds{1}_{(|\bar{\mathcal{E}}^m(\theta_i) / \sqrt{m}| > \epsilon)}.
\end{align*}
Pick $\delta > 0$.
\begin{align*}
    P \left( \frac{1}{n} \sum_{i=1}^n \mathds{1}_{(|\bar{\mathcal{E}}^m(\theta_i) / \sqrt{m}| > \epsilon)} > \delta \right) &\leq \frac{\mathbb{E}_{\theta, \mathbb{P}_{\theta}} \left[ \mathds{1}_{(|\bar{\mathcal{E}}^m(\theta_i) / \sqrt{m}| > \epsilon)} \right]}{\delta} \\
    &= \frac{P \left( |\bar{\mathcal{E}}^m(\theta_i) / \sqrt{m}| > \epsilon \right)}{\delta} \\
    &\leq \frac{P \left( (\bar{\mathcal{E}}^m(\theta_i))^4 > \epsilon^4 m^2 \right)}{\delta} \\
    &\leq \frac{\mathbb{E}_{\theta, \mathbb{P}_{\theta}} [(\bar{\mathcal{E}}^m(\theta_i))^4]}{\delta \epsilon^4 m^4}
\end{align*}
By assumption \ref{zhu-assumption}, $\mathbb{E}_{\theta, \mathbb{P}_{\theta}} [(\bar{\mathcal{E}}^m(\theta_i))^4] < \infty$. Then, by Borel-Cantelli lemma, $(*) \rightarrow 0$ w.p.1 as $n,m \rightarrow \infty$.
We have that
\begin{equation*}
    \limsup_{n,m \rightarrow \infty} \frac{1}{n} \sum_{i=1}^n | \mathds{1}_{(\hat{H}^m_i \geq \hat{v}^{n,m}_{\alpha})} - \mathds{1}_{(H_i \geq v_{\alpha})} | \leq 0.
\end{equation*}
Noting that $\frac{1}{n} \sum_{i=1}^n | \mathds{1}_{(\hat{H}^m_i \geq \hat{v}^{n,m}_{\alpha})} - \mathds{1}_{(H_i \geq v_{\alpha})} | \geq 0$ due to absolute value, the proof is complete.
\end{proof}



\begin{manualproposition}{3.12} \label{cvar-bias-proposition}
Under Assumptions \ref{zhu-assumption}, \ref{interchange-assumption}, \ref{supremum-assumption}, \ref{hong-assumption}, the bias $\mathbb{E}_{\theta, \mathbb{P}_{\theta}}[\psi^{n,m}_{\alpha}(x)] - c'_{\alpha}(x) \rightarrow 0$ as $n,m \rightarrow \infty$.
Moreover, if in addition $n = \Theta(m)$, then the bias is of the order $\mathcal{O}(n^{-1/2})$.
\end{manualproposition}

\begin{proof}

Note that $\mathbb{E}_{\theta} \left[ \frac{1}{(1-\alpha)n} \sum_{i=1}^n D_i \mathds{1}_{( H_i \geq v_{\alpha} )} \right] = c'_{\alpha}$. Therefore we will work with \linebreak $\mathbb{E}_{\theta} \left[ \frac{1}{(1-\alpha)n} \sum_{i=1}^n D_i \mathds{1}_{( H_i \geq v_{\alpha} )} \right]$ instead of $c'_{\alpha}$. 

\begin{align*}
    \mathbb{E}_{\theta, \mathbb{P}_{\theta}}[\psi^{n,m}_{\alpha}] - \mathbb{E}_{\theta} &\left[ \frac{1}{(1-\alpha)n} \sum_{i=1}^n D_i \mathds{1}_{( H_i \geq v_{\alpha} )} \right] 
    = \\
    &\mathbb{E}_{\theta, \mathbb{P}_{\theta}} \left[ \frac{1}{1-\alpha} \frac{1}{n} \sum_{i=1}^n \hat{D}^m_i \mathds{1}_{( \hat{H}^m_i \geq \hat{v}^{n,m}_{\alpha} )} - \frac{1}{1-\alpha} \frac{1}{n} \sum_{i=1}^n D_i \mathds{1}_{( H_i \geq v_{\alpha} )} \right] \\
    &\begin{aligned}
        = \frac{1}{1-\alpha} \mathbb{E}_{\theta, \mathbb{P}_{\theta}} &\left[ \frac{1}{n} \sum_{i=1}^n \hat{D}^m_i \mathds{1}_{( \hat{H}^m_i \geq \hat{v}^{n,m}_{\alpha} )} - \frac{1}{n} \sum_{i=1}^n D_i \mathds{1}_{( \hat{H}^m_i \geq \hat{v}^{n,m}_{\alpha} )} \right. \\
        &\left. + \frac{1}{n} \sum_{i=1}^n D_i \mathds{1}_{( \hat{H}^m_i \geq \hat{v}^{n,m}_{\alpha} )} - \frac{1}{n} \sum_{i=1}^n D_i \mathds{1}_{( H_i \geq v_{\alpha} )} \right] \\
    \end{aligned} \\
    &\begin{aligned}
        = \frac{1}{1-\alpha}& \left( \underbrace{ \mathbb{E}_{\theta, \mathbb{P}_{\theta}} \left[ \frac{1}{n} \sum_{i=1}^n (\hat{D}^m_i - D_i) \mathds{1}_{( \hat{H}^m_i \geq \hat{v}^{n,m}_{\alpha} )} \right] }_{(P1)} \right. \\
        &+\left. \underbrace{ \mathbb{E}_{\theta, \mathbb{P}_{\theta}} \left[ \frac{1}{n} \sum_{i=1}^n D_i ( \mathds{1}_{( \hat{H}^m_i \geq \hat{v}^{n,m}_{\alpha} )} - \mathds{1}_{( H_i \geq v_{\alpha} )} ) \right] }_{(P2)} \right)
    \end{aligned}
\end{align*}
We will show that both parts go to 0 w.p.1.
\begin{align*}
    |(P1)| &\leq \mathbb{E}_{\theta, \mathbb{P}_{\theta}} \left[ \frac{1}{n} \sum_{i=1}^n |\hat{D}^m_i - D_i| \mathds{1}_{( \hat{H}^m_i \geq \hat{v}^{n,m}_{\alpha} )} \right] \leq \mathbb{E}_{\theta, \mathbb{P}_{\theta}} \left[ \frac{1}{n} \sum_{i=1}^n |\hat{D}^m_i - D_i| \right] = \mathbb{E}_{\theta, \mathbb{P}_{\theta}}[ |\hat{D}^m (\theta) - D (\theta)| ]
\end{align*}
We have that for any $\theta$, $\mathbb{E}_{\mathbb{P}_{\theta}}[ |\hat{D}^m(\theta) - D(\theta)| ] \rightarrow 0$ w.p.1 by the strong consistency of $\hat{D}^m$. We will show how this carries on to $\mathbb{E}_{\theta, \mathbb{P}_{\theta}}[ |\hat{D}^m (\theta) - D (\theta)| ]$ using the assumption $\sup_{\theta} \mathbb{E}_{\mathbb{P}_{\theta}} [d(x, \xi)^2] < \infty$.

\begin{align*}
    \mathbb{E}_{\theta, \mathbb{P}_{\theta}}[ |\hat{D}^m (\theta) - D (\theta)| ] &\leq \sup_{\theta} \mathbb{E}_{\mathbb{P}_{\theta}}[ |\hat{D}^m (\theta) - D (\theta)| ] = \sup_{\theta} \mathbb{E}_{\mathbb{P}_{\theta}} \left[ \left| \frac{1}{m} \sum_{j=1}^m d(\xi_j) - D(\theta)  \right| \right] \\
    &\leq \sup_{\theta} \sqrt{\mathbb{E}_{\mathbb{P}_{\theta}} \left[ \left( \frac{1}{m} \sum_{j=1}^m d(\xi_j) - D(\theta)  \right)^2 \right]} = \sup_{\theta} \sqrt{\frac{\text{Var}_{\mathbb{P}_{\theta}}(d(\xi))}{m}} \\
    &\leq \frac{1}{\sqrt{m}} \sup_{\theta} \sqrt{\mathbb{E}_{\mathbb{P}_{\theta}} [d(\xi)^2]} \rightarrow 0 
\end{align*}
as $m \rightarrow \infty$ since $\sup_{\theta} \mathbb{E}_{\mathbb{P}_{\theta}} [d(x, \xi)^2] < \infty$. Therefore $(P1) \rightarrow 0$ as $n,m \rightarrow \infty$.

\begin{align*}
    |(P2)| &\leq \mathbb{E}_{\theta, \mathbb{P}_{\theta}} \left[ \frac{1}{n} \sum_{i=1}^n |D_i| | \mathds{1}_{( \hat{H}^m_i \geq \hat{v}^{n,m}_{\alpha} )} - \mathds{1}_{( H_i \geq v_{\alpha} )} | \right] \leq \mathbb{E}_{\theta, \mathbb{P}_{\theta}} \left[ \frac{1}{n} \sum_{i=1}^n \sup_{\theta}|D(\theta)| | \mathds{1}_{( \hat{H}^m_i \geq \hat{v}^{n,m}_{\alpha} )} - \mathds{1}_{( H_i \geq v_{\alpha} )} | \right] \\
    &= \sup_{\theta}|D (\theta)| \underbrace{ \mathbb{E}_{\theta, \mathbb{P}_{\theta}} \left[ \frac{1}{n} \sum_{i=1}^n | \mathds{1}_{( \hat{H}^m_i \geq \hat{v}^{n,m}_{\alpha} )} - \mathds{1}_{( H_i \geq v_{\alpha} )} | \right] }_{(P3)}
\end{align*}

Note that the term inside the expectation is bounded by $1$. Then, by Proposition \ref{indicator-proposition} and Dominated Convergence Theorem, $(P3) \rightarrow 0$. $\sup_{\theta}|D(\theta)| < \infty$ follows from the assumption \linebreak $\sup_{\theta} \mathbb{E}_{\mathbb{P}_{\theta}} [d(x, \xi)^2] < \infty$. Thus, we get $(P2) \rightarrow 0$ as $n,m \rightarrow \infty$. Therefore the bias converges to $0$ as $n,m \rightarrow \infty$.

To obtain the convergence rate, we use the same decomposition. Define $\mathcal{E}'(\xi(\theta)) = d(\xi(\theta)) - D(\theta)$ as the zero mean error term. Note that this term has a bounded variance by Assumption \ref{supremum-assumption}. For $(P1)$, it was shown that
\begin{equation*}
    |(P1)| \leq \mathbb{E}_{\theta, \mathbb{P}_{\theta}} \left[ \left| \hat{D}^m(\theta) - D(\theta) \right|  \right] = \mathbb{E}_{\theta, \mathbb{P}_{\theta}} \left[ \left| \frac{1}{m}\sum_{j=1}^m \mathcal{E}'(\xi_j(\theta)) \right|  \right].
\end{equation*}
Scaling with $\sqrt{m}$, we get
\begin{align*}
    \sqrt{m}\mathbb{E}_{\theta, \mathbb{P}_{\theta}} \left[ \left| \frac{1}{m}\sum_{j=1}^m \mathcal{E}'(\xi_j(\theta)) \right|  \right] &= \mathbb{E}_{\theta, \mathbb{P}_{\theta}} \left[ \left| \sqrt{m} \frac{1}{m}\sum_{j=1}^m \mathcal{E}'(\xi_j(\theta)) \right|  \right]
\end{align*}
where the term inside the absolute value converges to a mean zero normal random variable by CLT. It follows that $\sqrt{m}(P1) = \mathcal{O}(1)$ and $(P1) = \mathcal{O}(m^{-1/2})$.

We have the following lemma by \cite{Zhu2017input-uncertainty} that is useful when working with (P3). In the following $\breve{v}^m_{\alpha}$ is the $\alpha$ quantile of $\hat{H}^m$, i.e. $\hat{F}^m(\breve{v}^m_{\alpha}) = \alpha$ where $\hat{F}^m(\cdot)$ is the CDF of the noised response function and $\hat{f}^m(\cdot)$ is the PDF of it.
\begin{manuallemma}{A.1} Lemma B.4, \cite{Zhu2017input-uncertainty} - 
    Under Assumption \ref{zhu-assumption}, 
    \begin{equation*}
        \hat{v}^{n,m}_{\alpha} - \breve{v}^m_{\alpha} = \frac{1}{\hat{f}^m(\breve{v}^m_{\alpha})} \left( \alpha - \frac{1}{n} \sum_{i=1}^n \mathds{1}_{(\hat{H}^m(\theta_i) \leq \breve{v}^m_{\alpha})} \right) + A_n,
    \end{equation*}
    where $A_n = \mathcal{O}_{a.s.}(n^{-3/4} (\log n)^{3/4})$ and holds uniformly for all $m$.
\end{manuallemma}
\noindent Recall that 
\begin{equation*}
    (P2) = \mathbb{E}_{\theta, \mathbb{P}_{\theta}} \left[ \frac{1}{n} \sum_{i=1}^n D_i ( \mathds{1}_{( \hat{H}^m_i \geq \hat{v}^{n,m}_{\alpha} )} - \mathds{1}_{( H_i \geq v_{\alpha} )} ) \right] = \mathbb{E}_{\theta, \mathbb{P}_{\theta}} \left[ D_i ( \mathds{1}_{( \hat{H}^m_i \geq \hat{v}^{n,m}_{\alpha} )} - \mathds{1}_{( H_i \geq v_{\alpha} )} ) \right]. 
\end{equation*}
Since the difference of the indicators is $ \in \{-1, 0, 1\}$ and the gradient is bounded by Assumption \ref{supremum-assumption}, we can write
\begin{align*}
    - \sup_{\theta} |D(\theta)| \mathbb{E}_{\theta, \mathbb{P}_{\theta}} \left[ \mathds{1}_{( \hat{H}^m_i \geq \hat{v}^{n,m}_{\alpha} )} - \mathds{1}_{( H_i \geq v_{\alpha} )} \right]& \\
    \leq \mathbb{E}_{\theta, \mathbb{P}_{\theta}} \Big[ D_i ( & \mathds{1}_{( \hat{H}^m_i \geq \hat{v}^{n,m}_{\alpha} )} - \mathds{1}_{( H_i \geq v_{\alpha} )} ) \Big] \\
    &\leq \sup_{\theta} |D(\theta)| \underbrace{ \mathbb{E}_{\theta, \mathbb{P}_{\theta}} \left[ \mathds{1}_{( \hat{H}^m_i \geq \hat{v}^{n,m}_{\alpha} )} - \mathds{1}_{( H_i \geq v_{\alpha} )} \right] }_{(P4)}.
\end{align*}
It follows that $(P2) = \mathcal{O}((P4))$. In the following, we use a trick from Section 4.2 of \cite{Hong2009CVaR} to mitigate the dependency of $\hat{H}^m(\theta)$ and $\hat{v}^{n,m}_{\alpha}$. After applying the trick, we have the following where $\hat{v}^{(n-1), m}$ is calculated using $\theta_2, \theta_3, \ldots, \theta_n$ and is independent of $\hat{H}^m(\theta_1)$. Note that by definition $\mathbb{E}_{\theta, \mathbb{P}_{\theta}}[\mathds{1}_{( H_i \geq v_{\alpha} )}] = \alpha = \mathbb{E}_{\theta, \mathbb{P}_{\theta}} [\mathds{1}_{(\hat{H}^m(\theta_1) > \breve{v}^m_{\alpha})}]$.
\begin{align*}
    (P4) &= \mathbb{E}_{\theta, \mathbb{P}_{\theta}} \left[ \mathds{1}_{(\hat{H}^m(\theta_1) > \hat{v}^{(n-1), m}_{\alpha})} \right] - \alpha \\
    &= \mathbb{E}_{\theta, \mathbb{P}_{\theta}} \left[ \mathds{1}_{(\hat{H}^m(\theta_1) > \hat{v}^{(n-1), m}_{\alpha})} - \mathds{1}_{(\hat{H}^m(\theta_1) > \breve{v}^m_{\alpha})} \right] \\
    &= \mathbb{E}_{\hat{v}^{(n-1), m}_{\alpha}} \left[ \mathbb{E}_{\theta, \mathbb{P}_{\theta}} \left[ \mathds{1}_{(\hat{H}^m(\theta_1) > \hat{v}^{(n-1), m}_{\alpha})} - \mathds{1}_{(\hat{H}^m(\theta_1) > \breve{v}^m_{\alpha})} \mid \hat{v}^{(n-1), m}_{\alpha} \right] \right] \\
    &= \mathbb{E}_{\hat{v}^{(n-1), m}_{\alpha}} \left[ \hat{F}^m(\hat{v}^{(n-1), m}_{\alpha}) - \hat{F}^m(\breve{v}^m_{\alpha}) \right] \\
    &= \mathbb{E}_{\hat{v}^{(n-1), m}_{\alpha}} \left[ \hat{f}^m(\breve{v}^m_{\alpha}) (\hat{v}^{(n-1), m}_{\alpha} - \breve{v}^m_{\alpha}) + o(\hat{v}^{(n-1), m}_{\alpha} - \breve{v}^m_{\alpha}) \right].
\end{align*}
If we ignore the constant and the $o(\cdot)$ terms, we need to show that
\begin{equation}
    \sqrt{n} \mathbb{E}_{\hat{v}^{(n-1), m}_{\alpha}} \left[ (\hat{v}^{(n-1), m}_{\alpha} - \breve{v}^m_{\alpha}) \right] \rightarrow 0.
\end{equation}
It is seen from the Proof of Theorem 3.6 of \cite{Zhu2017input-uncertainty} that under Assumption \ref{zhu-assumption} $n = o(m^2)$ is a sufficient and necessary condition for 
\begin{equation}
    \lim_{n, m \rightarrow \infty} \sqrt{n} (\hat{v}^{(n-1), m}_{\alpha} - \breve{v}^m_{\alpha}) \Rightarrow \mathcal{N}(0, \sigma_v)
\end{equation}
where $\sigma_v = \frac{\alpha (1 - \alpha)}{f^2(v_{\alpha})}$. All that is left is to justify the interchange of the limit and expectation. For this, it suffices to show that $\sup_{n,m} \mathbb{E}_{\theta, \mathbb{P}_{\theta}}[n (\hat{v}^{(n-1), m}_{\alpha} - \breve{v}^m_{\alpha})^2] = \sup_{n,m} \mathbb{E}_{\theta, \mathbb{P}_{\theta}}[(n+1) (\hat{v}^{n, m}_{\alpha} - \breve{v}^m_{\alpha})^2] < \infty$, i.e. the sequence is uniformly integrable. 
\begin{align*}
    \mathbb{E}_{\theta, \mathbb{P}_{\theta}}[(n+1) (\hat{v}^{n, m}_{\alpha} - \breve{v}^m_{\alpha})^2] &= (n+1) \mathbb{E}_{\theta, \mathbb{P}_{\theta}} \left[ \left( \hat{f}^m(\breve{v}^m_{\alpha}) \left( \alpha - \frac{1}{n} \sum_{i=1}^n \mathds{1}_{(\hat{H}^m(\theta_i) \leq \breve{v}^m_{\alpha})} \right) + A_n \right)^2 \right] \\
    &= (n+1) (\hat{f}^m(\breve{v}^m_{\alpha}))^2 \mathbb{E}_{\theta, \mathbb{P}_{\theta}} \left[ \left( \alpha - \frac{1}{n} \sum_{i=1}^n \mathds{1}_{(\hat{H}^m(\theta_i) \leq \breve{v}^m_{\alpha})} \right)^2 \right] \\
    &\qquad \qquad + 2 (n+1) \hat{f}^m(\breve{v}^m_{\alpha}) \mathbb{E}_{\theta, \mathbb{P}_{\theta}} \left[  \alpha - \frac{1}{n} \sum_{i=1}^n \mathds{1}_{(\hat{H}^m(\theta_i) \leq \breve{v}^m_{\alpha})}  \right] A_n \\
    &\qquad \qquad + (n+1) A_n^2.
\end{align*}
Here, the expectation in the second term equals zero and the third term is bounded by definition. So, ignoring the constant terms, we only need to show that the following is bounded.
\begin{align*}
    (n+1) \mathbb{E}_{\theta, \mathbb{P}_{\theta}} \left[ \left( \alpha - \frac{1}{n} \sum_{i=1}^n \mathds{1}_{(\hat{H}^m(\theta_i) \leq \breve{v}^m_{\alpha})} \right)^2 \right] &= (n+1) \mathbb{E}_{\theta, \mathbb{P}_{\theta}} \left[ \left( \alpha - \frac{1}{n} \sum_{i=1}^n \mathds{1}_{(\hat{H}^m(\theta_i) \leq \breve{v}^m_{\alpha})} \right)^2 \right] \\
    &= (n+1) \mathbb{E}_{\theta, \mathbb{P}_{\theta}} \left[ \alpha^2 - 2 \alpha \frac{1}{n} \sum_{i=1}^n \mathds{1}_{(\hat{H}^m(\theta_i) \leq \breve{v}^m_{\alpha})} \right. \\
    & \qquad \qquad \left.+ \left( \frac{1}{n} \sum_{i=1}^n \mathds{1}_{(\hat{H}^m(\theta_i) \leq \breve{v}^m_{\alpha})} \right)^2 \right] \\
    &= (n+1) \left( \alpha^2 - 2 \alpha^2 + \frac{1}{n^2} (n \alpha + n (n-1) \alpha ^2) \right) \\
    &= \frac{n+1}{n} \alpha - \alpha^2
\end{align*}
where the results holds uniformly for all $m$ and $\sup_n \frac{n+1}{n} \alpha - \alpha^2 < \infty$. It follows that 
\begin{equation}
    \lim_{n, m \rightarrow \infty} \sqrt{n} \mathbb{E}_{\theta, \mathbb{P}_{\theta}} [\hat{v}^{(n-1), m}_{\alpha} - \breve{v}^m_{\alpha}] \rightarrow 0 \ \text{as} \ n \rightarrow \infty.
\end{equation}
Therefore, we get that $(P2) = o(n^{-1/2})$. Recalling the relation $n = \Theta(m)$ and putting the two together, we conclude that the bias is $\mathcal{O}(n^{-1/2})$.

\end{proof}



\begin{manualtheorem}{3.13} \label{consistency-theorem}
Under Assumptions \ref{zhu-assumption}, \ref{interchange-assumption}, \ref{supremum-assumption}, \ref{hong-assumption}, $\psi^{n,m}_{\alpha}(x)$ is a strongly consistent estimator of $c'_{\alpha}$. 
\end{manualtheorem}

\begin{proof}
Recall that $\psi^{n,m}_{\alpha} = \frac{1}{(1-\alpha) n} \sum_{i=1}^n \hat{D}^m_i \mathds{1}_{( \hat{H}^m_i \geq \hat{v}^{n,m}_{\alpha} )}$. We need to show that $\psi^{n,m}_{\alpha} - c'_{\alpha} \rightarrow 0$ w.p.1 as $n,m \rightarrow \infty$. 
\begin{align*}
\psi^{n,m}_{\alpha} - c'_{\alpha} &= \frac{1}{(1-\alpha)} \left[ \frac{1}{n} \sum_{i=1}^n \hat{D}^m_i \mathds{1}_{( \hat{H}^m_i \geq \hat{v}^{n,m}_{\alpha} )} - \mathbb{E}_{\theta} [D(\theta) \mathds{1}_{( H(\theta) \geq v_{\alpha} )} ] \right] \\
&= \frac{1}{(1-\alpha)} \left[ \underbrace{ \frac{1}{n} \sum_{i=1}^n \hat{D}^m_i \mathds{1}_{( \hat{H}^m_i \geq \hat{v}^{n,m}_{\alpha} )} - \frac{1}{n} \sum_{i=1}^n \hat{D}^m_i \mathds{1}_{( H_i \geq v_{\alpha} )} }_{(P1)} \right. \\
& \qquad \qquad \qquad  \left. + \underbrace{ \frac{1}{n} \sum_{i=1}^n \hat{D}^m_i \mathds{1}_{( H_i \geq v_{\alpha} )} - \mathbb{E}_{\theta} [D(\theta) \mathds{1}_{( H(\theta) \geq v_{\alpha} )} ] }_{(P2)} \right]
\end{align*}
We will show that each part individually goes to zero w.p.1.
\begin{align*}
(P2) &= \frac{1}{n} \sum_{i=1}^n \hat{D}^m_i \mathds{1}_{( H_i \geq v_{\alpha} )} - \mathbb{E}_{\theta} [D(\theta) \mathds{1}_{( H(\theta) \geq v_{\alpha} )} ] \\
&= \frac{1}{n} \sum_{i=1}^n \hat{D}^m_i \mathds{1}_{( H_i \geq v_{\alpha} )} - \frac{1}{n} \sum_{i=1}^n D_i \mathds{1}_{( H_i \geq v_{\alpha} )} + \frac{1}{n} \sum_{i=1}^n D_i \mathds{1}_{( H_i \geq v_{\alpha} )} - \mathbb{E}_{\theta} [D(\theta) \mathds{1}_{( H(\theta) \geq v_{\alpha} )} ]
\end{align*}
Here, the second part goes to zero with probability one by Strong Law of Large Numbers. For the first part, 
\begin{align*}
    \left| \frac{1}{n} \sum_{i=1}^n \hat{D}^m_i \mathds{1}_{( H_i \geq v_{\alpha} )} - \frac{1}{n} \sum_{i=1}^n D_i \mathds{1}_{( H_i \geq v_{\alpha} )} \right|
    &\leq \frac{1}{n} \sum_{i=1}^n |\hat{D}^m_i - D_i| \mathds{1}_{( H_i \geq v_{\alpha} )} \\
    &\leq \sup_{\theta}  |\hat{D}^m (\theta) - D(\theta)| \\
    &= \sup_{\theta} \left| \frac{1}{m} \sum_{j=1}^m d(\xi(\theta)) - D(\theta) \right|
\end{align*}
which is the sample average of a mean $0$ random variable with finite variance (by the assumption $\sup_{\theta} \mathbb{E}_{\mathbb{P}_{\theta}} [d(x, \xi)^2] < \infty$). Therefore, by SLLN, $(P2) \rightarrow 0$ w.p.1 as $n,m \rightarrow \infty$. 
We need to show that $(P1) = \frac{1}{n} \sum_{i=1}^n \hat{D}^m_i ( \mathds{1}_{( \hat{H}^m_i \geq \hat{v}^{n,m}_{\alpha} )} -  \mathds{1}_{( H_i \geq v_{\alpha} )} ) \rightarrow 0$ w.p.1 as $n,m \rightarrow \infty$. By Cauchy-Schwarz inequality, we have
\begin{equation*}
    \left| \frac{1}{n} \sum_{i=1}^n \hat{D}^m_i ( \mathds{1}_{( \hat{H}^m_i \geq \hat{v}^{n,m}_{\alpha} )} -  \mathds{1}_{( H_i \geq v_{\alpha} )} ) \right| 
    \leq \left[\frac{1}{n} \sum_{i=1}^n |\hat{D}^m_i|^2 \right]^{1/2} \left[ \frac{1}{n} \sum_{i=1}^n | \mathds{1}_{( \hat{H}^m_i \geq \hat{v}^{n,m}_{\alpha} )} -  \mathds{1}_{( H_i \geq v_{\alpha} )} |^{2} \right]^{1/2}
\end{equation*}
In order to show that the first term is finite, we will show that $\frac{1}{n} \sum_{i=1}^n (\hat{D}^m_i)^{2} \rightarrow \mathbb{E}_{\theta}[D(\theta)^{2}]$ w.p.1 which is bounded by the assumption $\sup_{\theta} \mathbb{E}_{\mathbb{P}_{\theta}} [d(x, \xi)^2] < \infty$.
\begin{align*}
\frac{1}{n} \sum_{i=1}^n (\hat{D}^m_i)^{2} - \mathbb{E}_{\theta}[D(\theta)^{2}] &= \frac{1}{n} \sum_{i=1}^n (\hat{D}^m_i)^{2} - \frac{1}{n} \sum_{i=1}^n (D_i)^{2} + \frac{1}{n} \sum_{i=1}^n (D_i)^{2} - \mathbb{E}_{\theta}[D(\theta)^{2}] 
\end{align*}
The second part goes to zero w.p.1 by SLLN. We need to show that $\frac{1}{n} \sum_{i=1}^n ((\hat{D}^m_i)^{2} - (D_i)^{2}) \rightarrow 0$ w.p.1. Define $\mathcal{E}'(\xi(\theta)) := d(\xi(\theta)) - D(\theta)$.
\begin{align*}
    \frac{1}{n} \sum_{i=1}^n ((\hat{D}^m_i)^{2} - (D_i)^{2}) & \leq \sup_{\theta} \left( \hat{D}^m (\theta)^{2} - D(\theta)^{2} \right) \\
    &= \sup_{\theta} \left[ \left( \frac{1}{m} \sum_{j=1}^m \mathcal{E}'(\xi_j(\theta)) \right)^2 - 2 D(\theta) \frac{1}{m} \sum_{j=1}^m \mathcal{E}'(\xi_j(\theta)) \right]
\end{align*}
Note that $\sup_{\theta} D(\theta) < \infty$ and $\mathcal{E}'(\xi_j(\theta))$ is a mean zero random variable with finite variance by the assumption $\sup_{\theta} \mathbb{E}_{\mathbb{P}_{\theta}} [d(x, \xi)^2] < \infty$. Therefore, both terms converge to zero w.p.1 and we get $\left[\frac{1}{n} \sum_{i=1}^n |\hat{D}^m_i|^2 \right]^{1/2} < \infty$.

To complete the proof, we need to show that $\left[ \frac{1}{n} \sum_{i=1}^n | \mathds{1}_{( \hat{H}^m_i \geq \hat{v}^{n,m}_{\alpha} )} -  \mathds{1}_{( H_i \geq v_{\alpha} )} |^{2} \right]^{1/2} \rightarrow 0$ w.p.1 as $n,m \rightarrow \infty$. This follows from Continuous Mapping Theorem and Proposition \ref{indicator-proposition}. Therefore, $\psi^{n,m}_{\alpha} \rightarrow c'_{\alpha}$ w.p.1 as $n,m \rightarrow \infty$. 
\end{proof}


\bibliographystyle{apacite} 
\bibliography{ref.bib}